\documentclass{article}
\usepackage{fullpage}
\usepackage{mathtools}

\usepackage{amssymb, amsmath, amsfonts, amsthm, graphics, mathrsfs, enumerate}
\usepackage{hyperref}
\hypersetup{
    colorlinks=true,
    linkcolor=blue,
    filecolor=magenta,      
    urlcolor=cyan,
    citecolor = {blue}
    }
\usepackage[hmargin=1 in, vmargin = 1 in]{geometry}
\usepackage{multicol} 

\usepackage{graphicx}
\usepackage{float}   
\restylefloat{table} 

\usepackage{algorithm} 
\usepackage[noend]{algpseudocode}
\newcommand{\Input}{\algorithmicrequire}

\newcommand {\norm} [1] { \lVert #1 \rVert}

\DeclareMathAlphabet{\mathpzc}{OT1}{pzc}{m}{it}

\newcommand{\Del}{\nabla}

\newcommand{\be}{\begin{equation}}
\newcommand{\ee}{\end{equation}}

\newcommand{\re}{\mbox{Re}}
\newcommand{\im}{\mbox{Im}}
\usepackage{bm}    
\usepackage[titletoc,title]{appendix}
\usepackage[export]{adjustbox}

\usepackage{graphicx}
\usepackage{epstopdf}
\usepackage{amsthm}
\usepackage{amsmath}
\usepackage{amssymb}
\usepackage{caption}
\usepackage{subcaption}
\usepackage[all]{xy}
\usepackage{verbatim}
\usepackage{mwe,subcaption,tikz}
\tikzset{boximg/.style={remember picture,red,thick,draw,inner sep=0pt,outer sep=0pt}}

\newtheorem{remark}{Remark}
\newcommand{\bigO}{O}  
\graphicspath{{fig/}}

\providecommand{\norm}[1]{\lVert#1\rVert}

\title{Solution of Stokes flow in complex nonsmooth 2D geometries
  via a linear-scaling high-order adaptive integral equation scheme}

\author{Bowei Wu\footnotemark[1]~, Hai Zhu\footnotemark[1]\,, Alex Barnett\footnotemark[2]\, and Shravan Veerapaneni\footnotemark[1]}

\begin{document}
\maketitle
\renewcommand{\thefootnote}{\fnsymbol{footnote}}
\footnotetext[1]{Department of Mathematics, University of Michigan, Ann Arbor, MI 48109}
\footnotetext[2]{Flatiron Institute, Simons Foundation, New York, NY 10022}

\begin{abstract} We present a fast, high-order accurate and adaptive boundary integral scheme for solving the Stokes equations in complex---possibly nonsmooth---geometries in two dimensions. The key ingredient is a set of panel quadrature rules capable of evaluating weakly-singular, nearly-singular and hyper-singular integrals to high accuracy. Near-singular integral evaluation, in particular, is done using an extension of the scheme developed in 
J.~Helsing and R.~Ojala, {\it J. Comput. Phys.} {\bf 227} (2008) 2899--2921. The boundary of the given geometry is ``panelized'' automatically to achieve user-prescribed precision. We show that this adaptive panel refinement procedure works well in practice even in the case of complex geometries with large number of corners. In one example, for instance, a model 2D vascular network with 378 corners required less than 200K discretization points to obtain a 9-digit solution accuracy. 
\end{abstract}
\maketitle


\section{Introduction} \label{sc:intro}
Since the pioneering work of Youngren and Acrivos \cite{Youngren75, Youngren76}, boundary integral equation (BIE) methods have been widely used for studying various particulate Stokes flow systems including drops, bubbles, capsules, vesicles, red blood cells and swimmers (e.g., see recent works \cite{sorgentone2018highly, contrast, bryngelson2018global} and references therein). BIE methods exploit the linearity of Stokes equations and offer several advantages such as reduction in dimensionality, exact satisfaction of far-field boundary conditions, and ease of handling moving geometries. Development of fast algorithms such as the fast multipole methods (FMMs) \cite{greengard1996integral, Ying04, wang2007parallel, anna08, gimbutas2015computational, malhotra2015pvfmm}, Ewald summation methods \cite{saintillan2005smooth, lindbo2010spectrally, kumar2012accelerated, af2017fast}, and variants tailored to Stokes equations \cite{sangani1996n, zinchenko2000efficient, wang2006algorithms, wang2018treecode} further extended their scope for solving problems in physically-realistic parameter regimes. Consequently, they are being used to investigate problems in a wide range of physical scales, from microhydrodynamics of isolated particles to large-scale flows generated by suspensions of particles (e.g., see \cite{petascale, nazockdast2017fast, Yan2019}). 

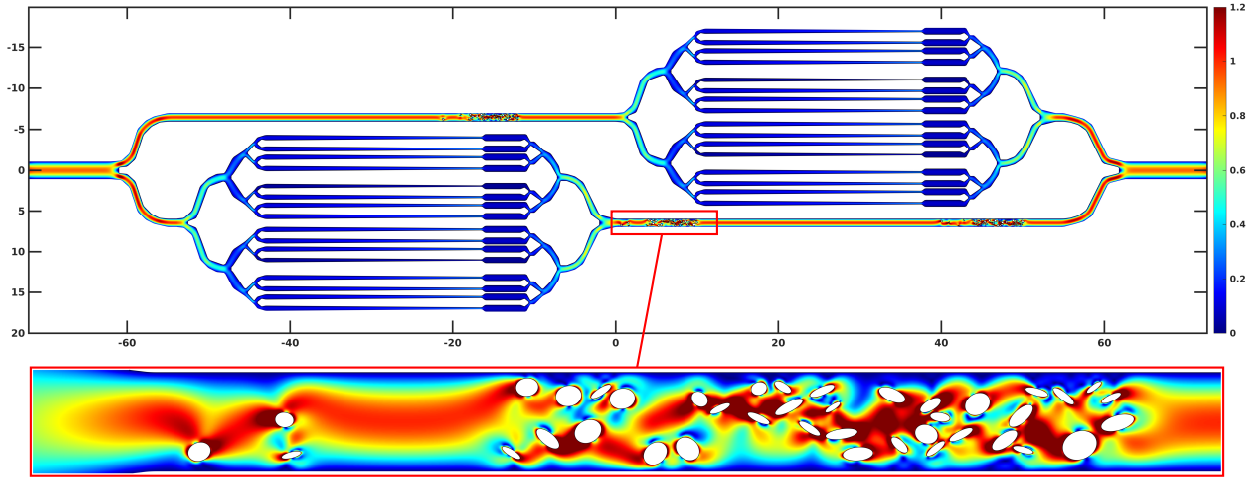
\begin{figure}[t]
  \centering
  \begin{subfigure}{\linewidth}
    \begin{tikzpicture}[boximg]
      \node[anchor=south west] (img) {\includegraphics[width=\linewidth]{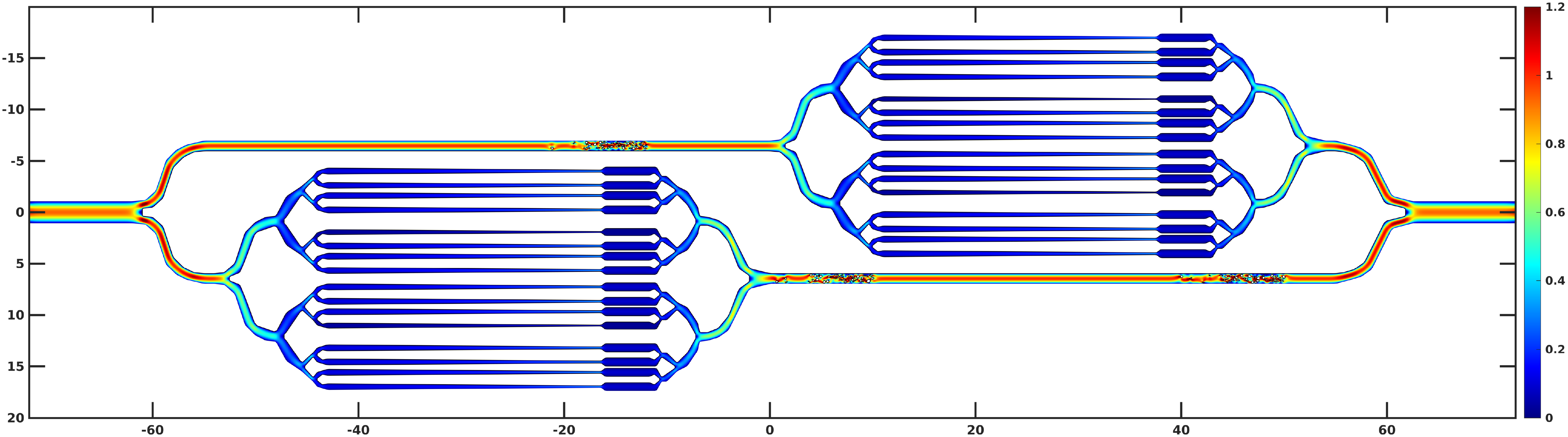}};
      \begin{scope}[x=(img.south east),y=(img.north west)]
        \node[draw,minimum height=0.3cm,minimum width=1.4cm] (B2) at (0.53,0.362) {};
      \end{scope}
    \end{tikzpicture}
  \end{subfigure}\\[0.5\baselineskip]
  \begin{subfigure}{1\linewidth}
   \centering
    \begin{tikzpicture}[boximg]
      \node (img2) {\includegraphics[width=0.96\linewidth]{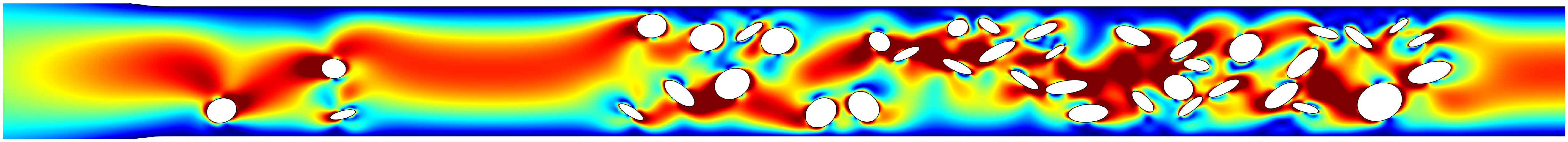}};
      \draw (img2.south west) rectangle (img2.north east);
    \end{tikzpicture}\hfill%

  \end{subfigure}
  \begin{tikzpicture}[overlay,boximg]
    \draw (B2) -- (img2);
  \end{tikzpicture}
  \caption{{\em Snapshot from a simulation of bacterial supension flow in a microfluidic chip geometry, which is inspired from the design proposed in \cite{hulme2007microfabricated}. A squirmer model \cite{ishikawa2007rheology} is used for modeling the bacteria, which treats them as rigid bodies with a prescribed slip at the fluid-structure interface. Thereby, we solve the Stokes equations with a no-slip boundary condition on the microfluidic chip geometry, a prescribed tangential velocity on the squirmer boundaries and an imposed parabolic flow profile at the inlet and outlet. We used $730,080$ discretization points for the chip boundary, resulting in $1,460,160$ degrees of freedom, and $128$ discretization points at each of the 120 squirmers. GMRES took about $10$ hours to reach a relative residual of $5.6\times 10^{-8}$, using an $8$-core $3.6$ GHz Intel Core i7 processor with 128 GB of RAM. Color indicates the magnitude of fluid velocity. 
  The estimated PDE relative $L_2-$norm error is $2\times 10^{-5}$.
  }}
 \label{fig:microfluidic_chip} 
\end{figure}

Despite their success, many numerical challenges still remain open for BIEs as applied to particulate Stokes flow problems. Prominent among them is the accurate handling of hydrodynamic interactions between surfaces that are almost in contact. For example, in dense suspension flows, the constituent particles often approach very close to each other or to the walls of the enclosing geometries (e.g., see Figure \ref{fig:microfluidic_chip}). To prevent artificial instabilities in this setting, numerical methods often require adaptive spatial discretizations, accurate nearly-singular integral evaluation schemes (note that while this issue is specific to BIEs, it manifests itself in other forms for other numerical methods) and accurate time-stepping schemes. The primary focus of this paper is addressing the first two issues for two-dimensional problems. 

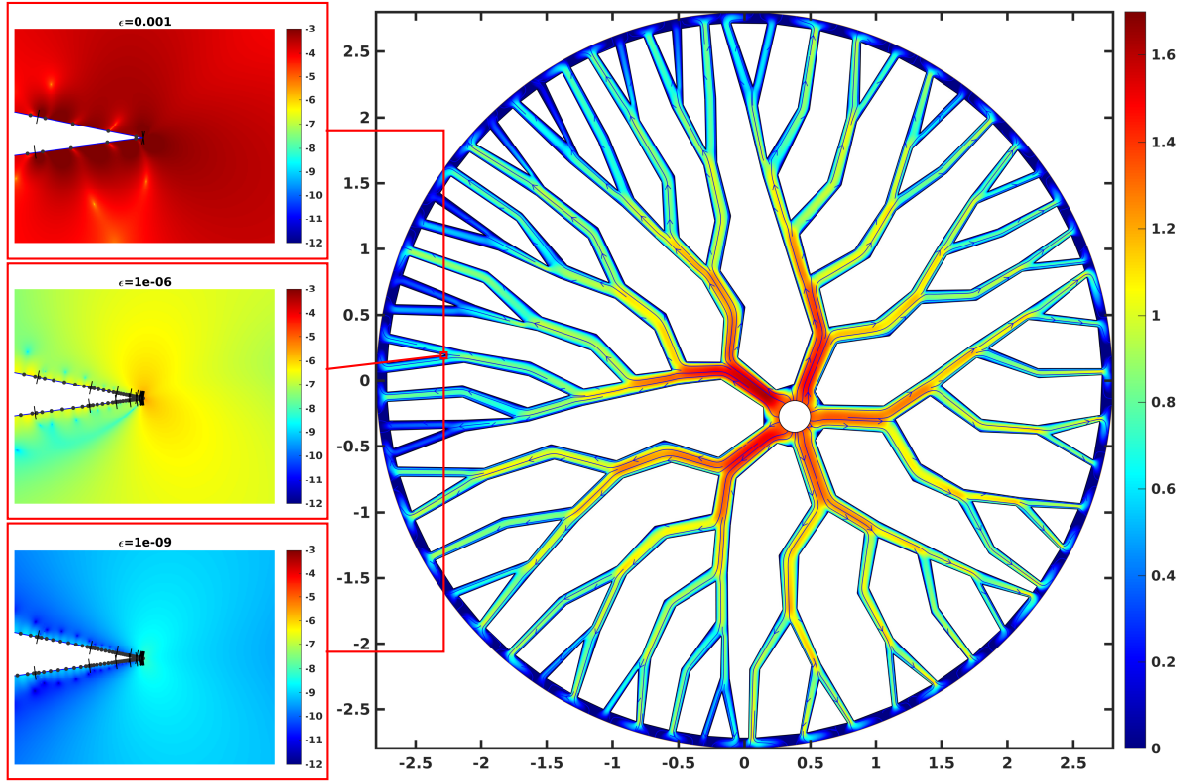
\begin{figure}[t]
  \centering
    \begin{subfigure}{0.3\linewidth}
    \begin{tikzpicture}[boximg]
      \node (img1) {\includegraphics[scale=0.34,right]{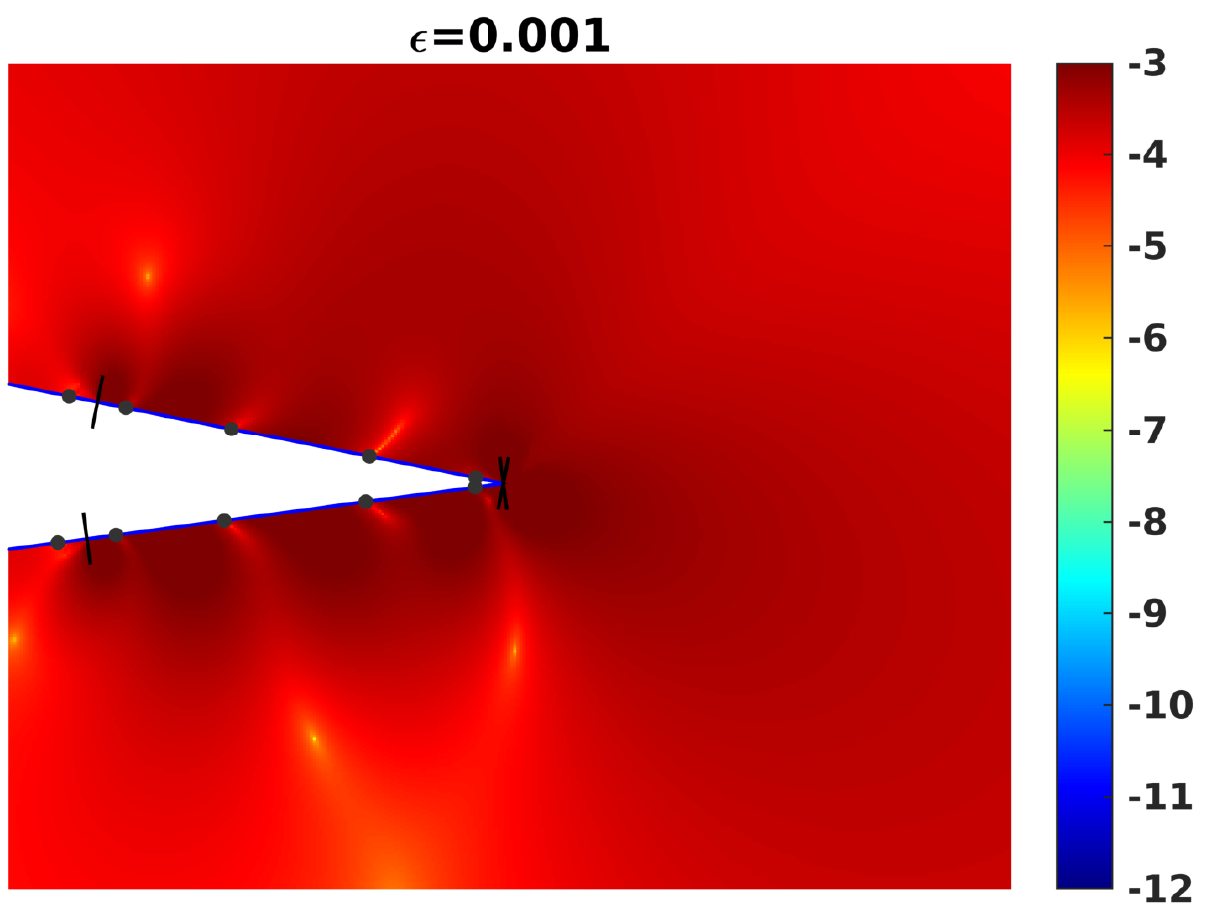}};
      \draw (-1.75,-1.7) rectangle (2.50,1.7);

    \end{tikzpicture}
    \\
    \begin{tikzpicture}[boximg]
      \node (img2) {\includegraphics[scale=0.34,right]{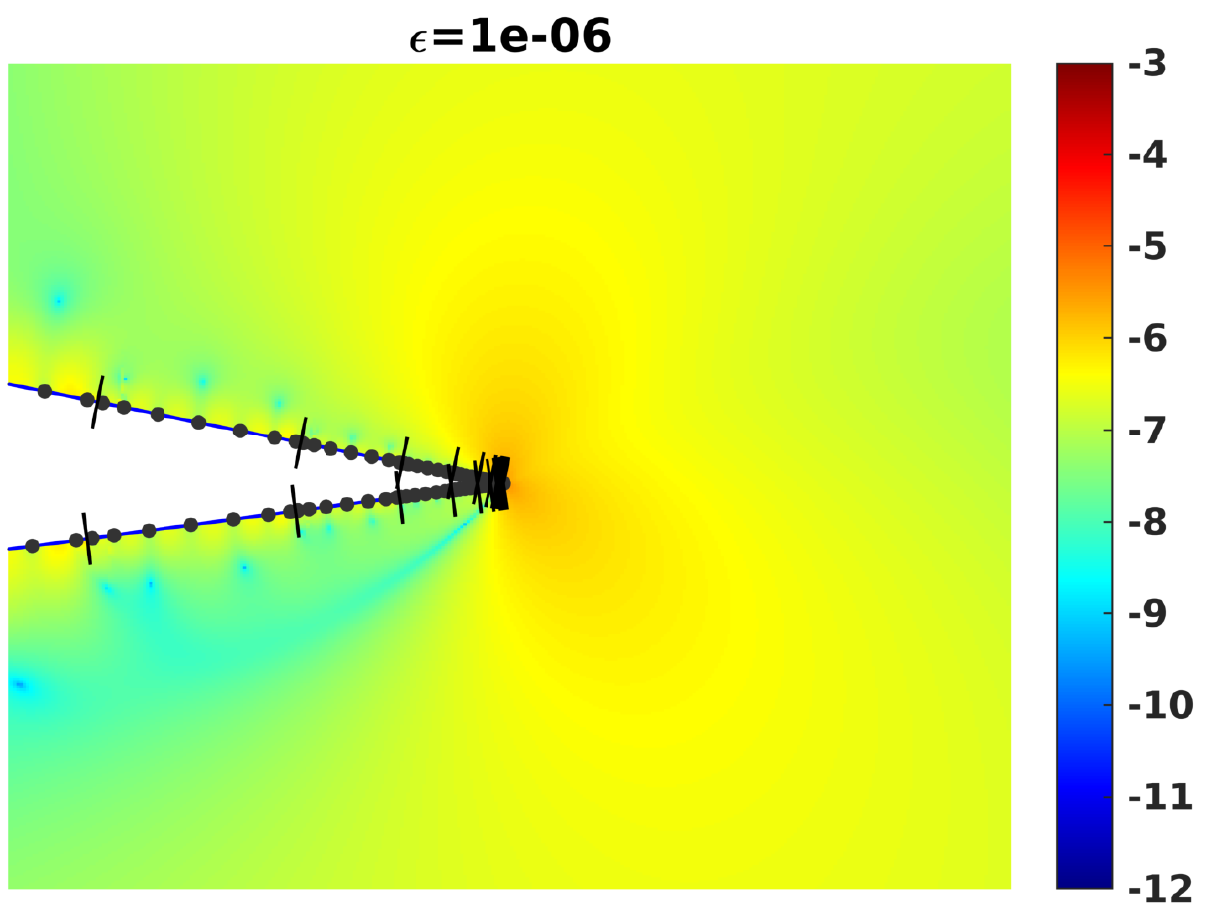}};
      \draw (-1.75,-1.7) rectangle (2.50,1.7);

    \end{tikzpicture}
    \\
    \begin{tikzpicture}[boximg]
      \node (img3) {\includegraphics[scale=0.34,right]{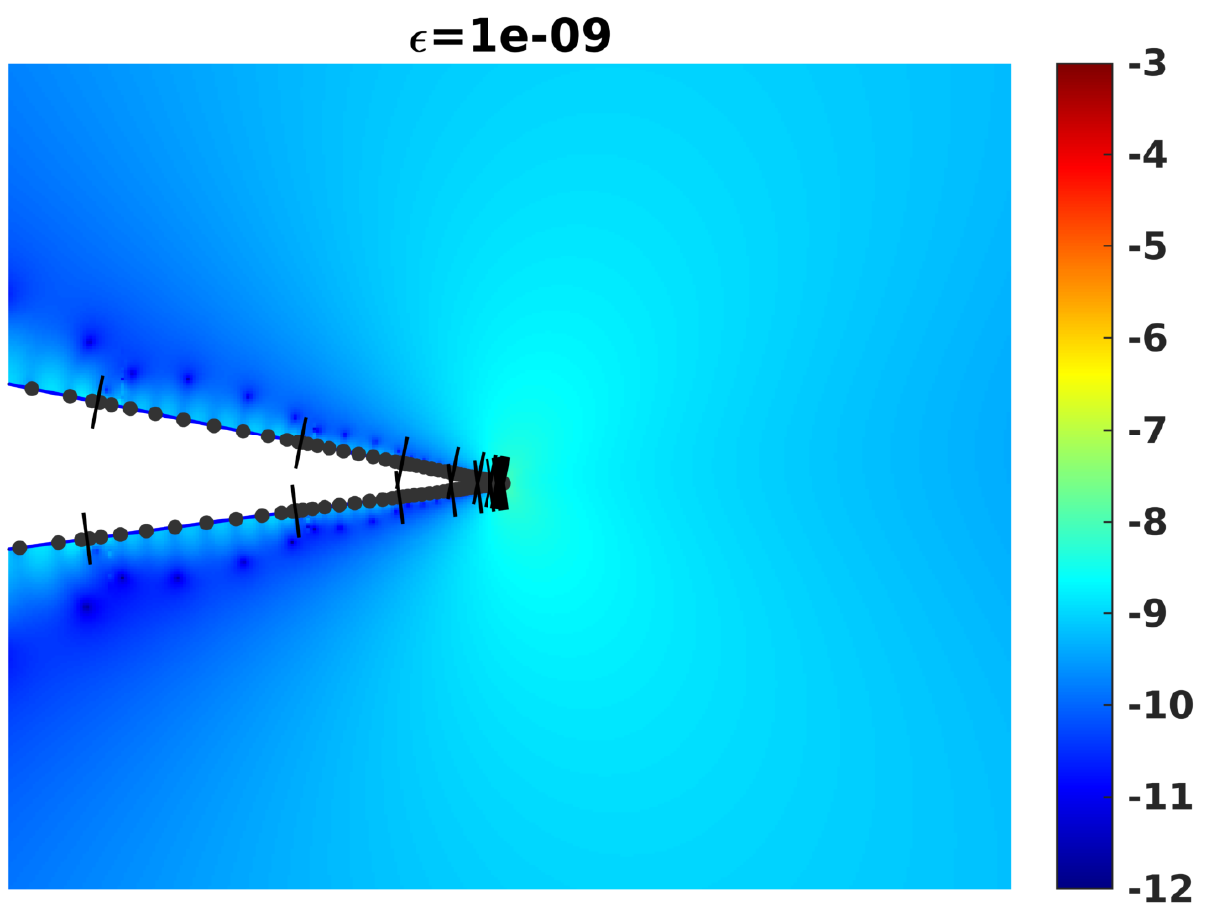}};
      \draw (-1.75,-1.7) rectangle (2.50,1.7);
    \end{tikzpicture}
    \label{fig:vasnet_panel}
  \end{subfigure}
  \begin{subfigure}{0.68\linewidth}
    \begin{tikzpicture}[boximg]
      \node[anchor=south west] (img) {\includegraphics[width=\linewidth]{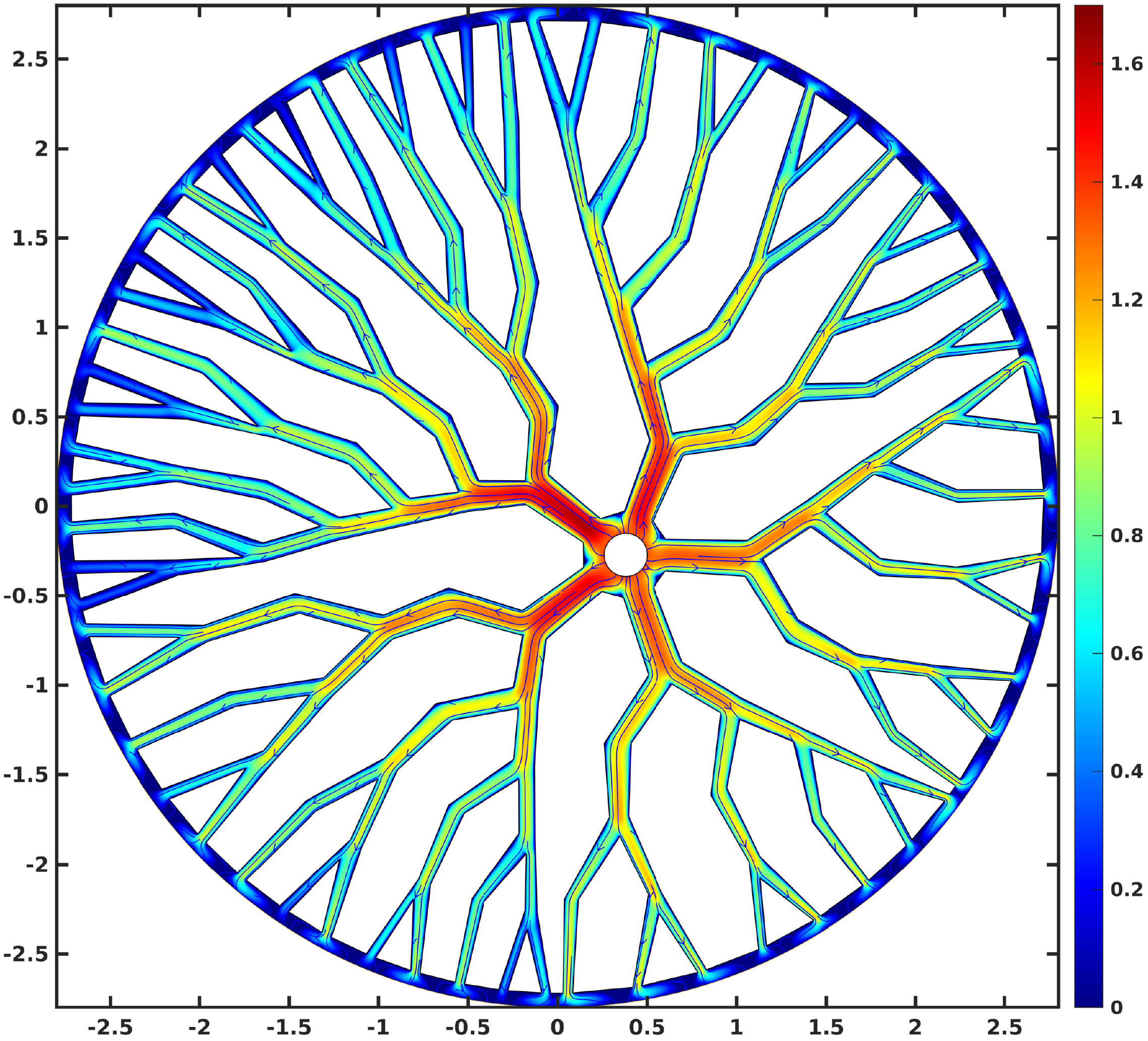}};
      \begin{scope}[x=(img.south east),y=(img.north west)]
        \node[draw,minimum height=0.08cm,minimum width=0.08cm] (B1) at (0.129,0.547) {};
      \end{scope}
    \end{tikzpicture}
    \label{fig:vasnet_flow}
  \end{subfigure}

  \begin{tikzpicture}[overlay,boximg]
    \draw (B1) |- (img1);
    \draw (B1) -- (img2);
    \draw (B1) |- (img3);
  \end{tikzpicture}
  \caption{{\em Solution of the Stokes equation in a nonsmooth circular vascular network with Dirichlet boundary condition. We apply no-slip boundary condition at all branch walls, and it is driven by a uniform flow from inner to outer circle. Color here indicates $\log$ of the magnitude of fluid velocity. We used automatically generated panels for both smooth boundaries and $378$ corners, resulting in $356,580$ degrees of freedom. GMRES took about $1$ hour to reach a relative residual of $7.61\times 10^{-11}$ on an $8$-core $4.0$ GHz Intel Core i7 desktop. The PDE solution has a relative $L_2-$norm error of $1\times 10^{-9}$. Three high-resolution $\log_{10}$ error plots that correspond to different user-requested tolerance $\epsilon$ near the same reentrant corner are shown on the left; here the short
  normal lines show panel endpoints, and the black dots quadrature nodes.}}
  \label{fig:vasnet}
\end{figure}

In this work, we introduce specialized panel quadrature schemes that can accurately evaluate layer potentials defined on a smooth open curve, and for target points arbitrarily close to, or on, the curve. This helps the efficient tackling of dense particulate flows constrained in large multiply-connected domains such as in Figure \ref{fig:microfluidic_chip}. In addition, we formulate a set of rules for refining (or coarsening) the panels used to represent the boundaries, so that a user-specified error tolerance can be achieved automatically. One of the advantages of this adaptive panel refinement procedure is that it can  handle geometries with corners (as in Figure \ref{fig:vasnet}) or nearly self-touching geometries, using the same quadrature schemes for both weakly- and nearly-singular integrals.

Our work is closely related to two recent efforts, that of Barnett et al.\ \cite{lsc2d} and Ojala--Tornberg \cite{ojala2015accurate}, both of which, in turn, are extensions to Stokes potentials of the Laplace work of Helsing--Ojala \cite{helsing-ojala}. In \cite{lsc2d}, three of the co-authors developed high-order global quadrature schemes based on the periodic trapezoidal rule (PTR) for the evaluation of Laplace and Stokes layer potentials defined on a smooth closed planar curve, which achieved spectral accuracy for particles arbitrarily close to each other. Nevertheless, a key limitation of these close evaluation schemes is that they only work for closed curves and have uniform resolution on any part of a curve. Thereby, while they are well-suited for moving rigid and deformable particles immersed in Stokes flow, new machinery is needed for the constraining complex geometries such as those in Figures \ref{fig:microfluidic_chip} and \ref{fig:vasnet}. 

On the other hand, Ojala--Tornberg \cite{ojala2015accurate} developed a panel quadrature scheme as done in this work. The main distinction is that \cite{ojala2015accurate} uses the BIE framework of Kropinski \cite{kropinski00}, converting Stokes equations into a biharmonic equation, then tailors it to a specific application (droplet hydrodynamics), whereas our scheme directly tackles all the Stokes layer potentials using physical variables (single-layer, double-layer and their associated pressure and traction kernels). Thereby, it can be integrated with existing BIE methods developed for various physical problems (colloids, drops, vesicles, squirmers and other suspensions) and flow conditions (imposed flow, pressure-, electrically- or magnetically-driven flows, etc) more naturally. 
We note that resolving nearly-singular integrals is an active area of research owing to its importance in solving several other linear elliptic partial differential equations (PDEs), such as Laplace, Helmholtz and biharmonic equations, via BIEs; other recent works that considered two-dimensional problems include \cite{klockner2013quadrature, barnett2014evaluation, carvalho2018asymptotic, rahimian2018ubiquitous, klinteberg2017adaptive, perez2019harmonic}.  

The predominant class of algorithms for adaptive meshing found in the Stokes BIE literature are the {\em reparameterization} schemes (also called the {\em resampling} techniques) that are dedicated to resolving boundary mesh distortions arising in deformable particle flow simulations (see \cite{ves3d} for a review on this topic). The primary focus of this work, on the other hand, is to determine an 
optimal distribution of boundary panels on a given domain so that a user-prescribed tolerance is met when computing the solution. Prior work in this area has mostly been restricted to low-order boundary element methods (BEMs); see \cite{kita2001error} for a review.
In the last 20 years, high-order $hp$-variants of BEM have been tested
for Laplace and Helmholtz problems on polygons
(for example \cite{heuer96,Bantle15}).
However, we are not aware of any Nystr\"om $hp$-BIEs for the Stokes equations capable of handling arbitrary complex geometries. The recent research of Rachh--Serkh \cite{rachh2017solution} exploits a power-law basis resulting from analysis of the wedge problem to solve Stokes corner problems via a BIE.
Here, we take a more pedestrian approach to handling corner singularities via the use of graded meshes. A prototypical example is shown in Figure \ref{fig:vasnet}. 

The paper is organized as follows. 
In Section \ref{sc:problem}, we define the boundary integral operators and state the integral equation formulations for the viscous flow problems considered here. In Section \ref{sc:panel}, we present the panel quadrature rules for evaluating the Stokes layer potentials. We summarize our algorithm for adaptively discretizing a given geometry in Section \ref{sc:amr}. We consider several test cases and present results on the performance of our algorithms in Section \ref{sc:results}, followed by conclusions in Section \ref{sc:conclusions}.  

\section{Mathematical preliminaries} \label{sc:problem}

In this section, we first state the PDE formulation for the example shown in Figure \ref{fig:vasnet}, reformulate it as a BIE and then discuss the evaluation of the resulting boundary integral operators. 
\subsection{Boundary value problem and integral equation formulation} 
\label{sc:formulation}
The fluid domain $\Omega$ in Figure \ref{fig:vasnet} is a multiply-connected region bounded by $N_\Gamma$ closed curves, $\{\Gamma_i, i = 1, \ldots, N_\Gamma\}$. Without loss of generality, let $\Gamma_1$ be the all-enclosing boundary---i.e., the outer circle in Figure \ref{fig:vasnet}---and let $\Gamma = \cup_i \Gamma_i$. Denoting the fluid viscosity by $\mu$, the velocity by $u$ and the pressure by $p$, the governing boundary-value problem, in the vanishing Reynolds number limit, is
\begin{subequations}
\begin{align}
-\mu \Delta u + \nabla p = 0 \quad\text{and}\quad \nabla \cdot u = 0 & \quad\text{in}\quad \Omega, \\
u = g &\quad\text{on}\quad \Gamma.
\end{align}
\label{PDE}
\end{subequations}
The Dirichlet data $g$ must satisfy the consistency condition $\int_\Gamma g \cdot n ds = 0$, where $n$ is the normal to $\Gamma$. For example, in Figure \ref{fig:vasnet}, the flow is driven by an outward flow condition at the inner circle and an inward flow condition at the outer circle. Their magnitude is chosen such that the consistency condition is respected. On the rest of the curves, a no-slip ($g=\bm{0}$) boundary condition is enforced. 

While there are many approaches for reformulating the Dirichlet problem (\ref{PDE}) as a BIE \cite{pozrikidis1992boundary,hsiao-wendland}, for simplicity, we use an indirect, combined-field BIE formulation that leads to a well-conditioned and non-rank-deficient linear system. We make the following {\em ansatz}:   
\be
u(x) = \sum_{i=1}^{N_\Gamma} \, ( \mathbf{S}_i \sigma) (x) + (\mathbf{D}_i \sigma )(x)
:= (\mathbf{S}\sigma)(x) + (\mathbf{D}\sigma)(x)~,
\qquad x \in \Omega,
\label{urep}
\ee
where $\sigma$ is an unknown vector ``density'' function to be determined, $\mathbf{S}_i$ is the velocity Stokes single layer potential (SLP) and $\mathbf{D}_i$ is the double layer potential (DLP), defined by  
\begin{align}
  (\mathbf{S}_i\sigma)(x) &= \frac{1}{4\pi\mu}\int_{\Gamma_i}\left(I \log\frac{1}{\rho} + \dfrac{r\otimes r}{\rho^2}\right)\sigma(y)\,ds_y, &x\in\Omega,
  \label{sto_slp}
  \\
  (\mathbf{D}_i\sigma)(x) &= \frac{1}{\pi}\int_{\Gamma_i}\left(\dfrac{r\cdot n_y}{\rho^2} \dfrac{r\otimes r}{\rho^2}\right)\sigma(y)\,ds_y, &x\in\Omega.
  \label{sto_dlp}
\end{align}
Here, $I$ is the 2-by-2 identity tensor, $r = x - y$, $\rho = \norm{r}_2$, and $ds_y$ is the arc length element on $\Gamma$.
As \eqref{urep} indicates, we use $\mathbf{S}$ or $\mathbf{D}$
to denote the sum of layer potentials due to all source curves.
When combined with their associated pressure kernels (given in the Appendix),
the convolution kernels above are fundamental solutions to the Stokes equations (\ref{PDE}a); therefore, all that is left to ensure that the {\em ansatz} (\ref{urep}) solves the boundary-value problem
is to enforce the velocity boundary condition (\ref{PDE}b).

We introduce the operator block notation $\mathbb{S}_{ij}\sigma = (\mathbf{S}_j \sigma) (\Gamma_i)$, and similarly $\mathbb{D}_{ij}$ for the double layer potential, i.e., the layer velocity potential restricted to the target curve $\Gamma_i$, using the source curve $\Gamma_j$.
By taking the limit as $x$ approaches $\Gamma$ (from the interior for $\Gamma_1$, and exterior for the rest of the curves),
the standard jump relations for the single and double layer potentials \cite{pozrikidis1992boundary,hsiao-wendland} give the
interior case of the following $N_\Gamma \times N_\Gamma$ BIE system:
\be
\left(
\begin{bmatrix}
  I/2 & 0 & \dots \\
  0 & I/2 & \dots \\
  \vdots &\vdots &\ddots
\end{bmatrix}
+
\begin{bmatrix}
  \eta(\mathbb{S}_{11}+\mathbb{D}_{11}) & \mathbb{S}_{12} + \mathbb{D}_{12} & \dots \\
  \eta(\mathbb{S}_{21}+\mathbb{D}_{21}) & \mathbb{S}_{22} + \mathbb{D}_{22} & \dots \\
  \vdots & \vdots & \ddots
\end{bmatrix}
\right)
\begin{bmatrix}
  \eta\sigma(\Gamma_1)\\
  \sigma(\Gamma_2)\\
  \vdots
\end{bmatrix}
=
\begin{bmatrix}
  f(\Gamma_1)\\
  f(\Gamma_2)\\
  \vdots
\end{bmatrix},
\quad
\eta = \left\{\begin{array}{ll}-1, & \mbox{interior case,}\\
+1, & \mbox{exterior case.}
\end{array}\right.
\label{bie}
\ee
The admixture of SLP and DLP in \eqref{urep}
insures that \eqref{bie} is of Fredholm second kind, and
that interior curves $\Gamma_i$, $i=2,3,\dots$ introduce no null-spaces
\cite[Thm.~2.1]{hebeker} \cite[p.128]{pozrikidis1992boundary}.
Note that the negation of the density for the outer curve
(first block column)
results in all positive identity blocks. Once we obtain the unknown density function $\sigma$ on $\Gamma$, by solving \eqref{bie}, we can evaluate the velocity at any point in the fluid domain by using \ref{urep}. 

We also consider the {\em exterior} boundary-value problem,
which can be thought of taking the above outer curve $\Gamma_1$ to infinity,
removing it from the problem.
The fluid domain $\Omega$
becomes the entire plane minus the closed interiors of the
other curves; these curves we may relabel as $\Gamma = \cup_{i=1}^{N_\Gamma} \Gamma_i$.
There is also now a given imposed background flow $u_\infty(x)$,
which in applications is commonly a uniform, shear, or extensional flow.
The decay condition
$u(x)-u_\infty(x) = \bm\Sigma \log \|x\| + O(1)$,
for some $\bm\Sigma\in\mathbb{R}^2$,
must be appended to \eqref{bie} to give a well-posed BVP \cite{hsiao-wendland}.
The representation \eqref{urep} of the physical velocity is now
\be
u(x) = u_\infty(x) + (\mathbf{D}\sigma+\mathbf{S}\sigma)(x)~.
\label{urepext}
\ee
The resulting BIE is given as the second case of \eqref{bie}.
In the simplest case of physical no-slip boundary conditions,
the data in \eqref{bie} is now $f=-u_\infty|_\Gamma$,
which one may check cancels the velocity on all boundaries.

The close-evaluation schemes of Helsing--Ojala \cite{helsing-ojala}
enable efficient and accurate evaluation of certain complex contour
integrals arbitrarily close to the boundary, in the case where high order panel-based discretization is used. Thus, a route to evaluate the needed Stokes layer potentials close to the boundary, as in \eqref{bie} and \eqref{urep}, is to express them in terms of Laplace potentials,
which are in turn expressed via contour integrals.
The rest of this section presents these formulae.

\subsection{Fundamental contour integrals}
\label{s:cont}
Let $\tau$ be a given, possibly complex, scalar density function on $\Gamma$.
We associate $\mathbb{R}^2$ with the complex plane $\mathbb{C}$.
Let $n_y$ be the outward-pointing unit vector at $y\in\Gamma$,
expressed as a complex number.
Given a target point $x\in\Omega$,
we define the potentials
\begin{align}
  I_L &= (I_L\tau)(x):=\int_{\Gamma}\log|x-y| \tau(y) |dy|
  = \int_{\Gamma}\log|x-y| \frac{\tau(y)}{in_y} dy
  ~, &\mbox{(real logarithmic)}
\label{IL}
\\
I_C &= (I_C\tau)(x):= \int_{\Gamma}\frac{\tau(y)}{y-x}dy~, &\mbox{(Cauchy)}
\label{IC}
\\
I_H&= (I_H\tau)(x) := \int_{\Gamma}\frac{\tau(y)}{(y-x)^2}dy~,  &\mbox{(Hadamard)}
\label{IH}
\\
I_S&= (I_S\tau)(x) := \int_{\Gamma}\frac{\tau(y)}{(y-x)^3}dy~.  &\mbox{(supersingular)}
\label{IS}
\end{align}
Note that, as functions of target point $x$, $I_C$, $I_H$ and $I_S$
are holomorphic functions in $\Omega$,
that $(d/dx)(I_C\tau)(x) = (I_H\tau)(x)$,
and $(d/dx)(I_H\tau)(x) = 2(I_S\tau)(x)$.
In contrast, $I_L$ is not generally holomorphic; yet, for $\tau$ real, $I_L$ is the real part of a holomorphic function.

\subsection{Laplace layer potentials}
Now, let $\tau$ be a real, scalar, density function on $\Gamma$.
The Laplace single- and double-layer potentials are defined, respectively, by 
\begin{equation}
  (\mathcal{S}\tau)(x):=\frac{1}{2\pi}\int_{\Gamma}\left(\log\frac{1}{\rho}\right)\tau(y)ds_y,\qquad x\in\Omega,
  \label{lap_slp}
\end{equation}
\begin{equation}
  (\mathcal{D}\tau)(x):=\frac{1}{2\pi}\int_{\Gamma}\left(\dfrac{\partial}{\partial n_y}\log\frac{1}{\rho}\right)\tau(y)ds_y = \frac{1}{2\pi}\int_{\Gamma}\left(\frac{r\cdot n_y}{\rho^2}\right)\tau(y)ds_y,\qquad x\in\Omega.
  \label{lap_dlp}
\end{equation}
In terms of the contour integrals \eqref{IL} and \eqref{IC},
using $ds_y = |dy| = dy/in_y$,
and the restriction to $\tau$ real,
we can rewrite the SLP and DLP as
\begin{equation}
  (\mathcal{S}\tau)(x) =\frac{-1}{2\pi} (I_L \tau)(x) \quad\text{and}\quad 
  (\mathcal{D}\tau)(x) = \re \, \frac{i}{2\pi} (I_C \tau)(x), \quad
  x\in\Omega.
\end{equation}

We will also need the gradients and Hessians of Laplace layer potentials.
The gradient of the SLP has components
\begin{equation}
  \dfrac{\partial}{\partial x_1}(\mathcal{S}\tau)(x) =
  \dfrac{1}{2\pi}\re\, (I_C[\tau/in_y])(x),
\qquad
\dfrac{\partial}{\partial x_2}(\mathcal{S}\tau)(x) =
\dfrac{-1}{2\pi}\im\, (I_C[\tau/in_y])(x)~.
\end{equation}
The gradient of the DLP requires the Hadamard kernel, and has components
\begin{equation}
  \dfrac{\partial}{\partial x_1}(\mathcal{D}\tau)(x) = \frac{-1}{2\pi}\im\, (I_H \tau)(x),
  \qquad
  \dfrac{\partial}{\partial x_2}(\mathcal{D}\tau)(x) = \frac{-1}{2\pi}\re\, (I_H \tau)(x)~.
  \label{gradD}
\end{equation}

The Hessians, which also require $I_S$, are discussed in Appendix~\ref{a:trac}.

\subsection{Stokes velocity layer potentials}

As in \cite{lsc2d}, we rewrite the
Stokes SLP in terms of Laplace potentials, and the DLP in terms
of Laplace and Cauchy potentials.
Using the identity
\[ \dfrac{ r\otimes r}{\rho^2}\sigma = \dfrac{ r}{\rho^2}( r\cdot\sigma) = ( r\cdot\sigma)\Del_{x}\log\rho, \]
we can rewrite the Stokes SLP in terms of the Laplace SLP \eqref{lap_slp}
as
\begin{equation}
\begin{aligned}
(\mathbf{S} \sigma)(x) \;=\;& {1\over4\pi\mu}\Big\{\int_{\Gamma}\left(\log\frac{1}{\rho}\right) \, \sigma ds_{y} + \Del\int_{\Gamma}\left(\log\frac{1}{\rho}\right) \, (y\cdot\sigma) ds_{y}\\
 & - x_{1}\Del\int_{\Gamma}\left(\log\frac{1}{\rho}\right) \, \sigma_{1}ds_{y}  - x_{2}\Del\int_{\Gamma}\left(\log\frac{1}{\rho}\right) \, \sigma_{2}ds_{y}\Big\}
~,
\end{aligned}
 \label{slp_stokes_via_lps}
\end{equation} 
where $\nabla = \nabla_x$ is assumed from now on.
Therefore, three Laplace potentials (and their first derivatives)
with real scalar density functions $y\cdot\sigma$, $\sigma_1$, and $\sigma_2$
need to be computed to evaluate the Stokes SLP. Similarly, using the identity
\begin{equation}
\Del\left ( \dfrac{ r\cdot n_{y}}{\rho^2} \right )
=\dfrac{ n_{y}}{\rho^2} - ( r\cdot n_{y})\dfrac{2 r}{\rho^4}
~,\label{identity_for_dlp_split}
\end{equation}
the Stokes DLP can be written as
\begin{equation}
\begin{aligned}
(\mathbf{D} \sigma)(x)\;
=\;& {1\over2\pi}\int_{\Gamma}{ n_{y}\over\rho^2}( r\cdot\sigma)ds_{y}  + {1\over2\pi}\Del\int_{\Gamma}{ r\cdot n_{y}\over\rho^2}( y\cdot\sigma)ds_{y}\\
 & - {1\over2\pi}x_{1}\Del\int_{\Gamma}{ r\cdot n_{y}\over\rho^2}\sigma_{1}ds_{y}  - {1\over2\pi}x_{2}\Del\int_{\Gamma}{r\cdot n_{y}\over\rho^2}\sigma_{2}ds_{y}
~.
\end{aligned}
\label{dlp_stokes_via_lps}
\end{equation}
While the last three terms are gradients of Laplace DLPs, the
first term requires a Cauchy integral.
More concisely, we can write \eqref{slp_stokes_via_lps} and \eqref{dlp_stokes_via_lps} as
\begin{align}
(\mathbf{S}\sigma)(x) &= \frac{1}{2\mu}\Big(\, (\mathcal{S}\sigma_1, \mathcal{S}\sigma_2)  + \nabla\mathcal{S}[y\cdot \sigma] - x_1\nabla \mathcal{S}\sigma_1 - x_2\nabla \mathcal{S}\sigma_2  \Big)(x), & x\in\Omega,   \label{eq:split_S}\\
(\mathbf{D}\sigma)(x) &= \Big(\, \frac{1}{2\pi} \text{Re}\,\big(I_C(\tau_1), I_C(\tau_2)\big) + \nabla \mathcal{D}[y\cdot \sigma] - x_1\nabla \mathcal{D}\sigma_1 - x_2\nabla \mathcal{D}\sigma_2\Big) (x), & x\in\Omega, \label{eq:split_D}
\end{align}
where, as above, a pair $(\cdot,\cdot)$ indicates two vector components,
and a short calculation verifies that
the complex scalar density functions $\tau_1$ and $\tau_2$
\begin{equation}
\tau_1 = (\sigma_1+i\sigma_2)\frac{\text{Re}\,n_y}{n_y},\qquad \tau_2 = (\sigma_1+i\sigma_2)\frac{\text{Im}\,n_y}{n_y},
\end{equation}
where $n_y$ is interpreted as a complex number,
makes the identity hold (i.e.\ the first term in
\eqref{eq:split_D} equals the first term in \eqref{dlp_stokes_via_lps}).

In summary, the procedure discussed in this section enables us to express the velocity field anywhere in the fluid domain, represented by \eqref{urep}, in terms of fundamental contour integrals.
Similar formulae exist for the fluid pressure and hydrodynamic stresses,
which are commonly needed in several applications;
we present these in Appendix~\ref{a:trac}.


\section{Nystr\"om discretization and evaluation of layer potentials}\label{sc:panel}

\subsection{Overview: discretization and the plain Nystr\"om formula}
\label{nystrom}

Firstly, we need to specify a numerical approximation of the density function $\sigma$.
For simplicity,
consider the case of an exterior BVP on a single closed curve
$\Gamma$ parameterized by $Z: [0,2\pi) \to \mathbb{R}^2$,
such that $\Gamma = Z([0,2\pi))$.
The Stokes BIE \eqref{bie} is then
\be
(I/2 + \mathbb{K})\sigma = f, \qquad \mathbb{K} := \mathbb{D}+\mathbb{S}~,
\label{biesimple}
\ee
where $\sigma$ and $f$ are 2-component vector functions on $\Gamma$.

Given the user requested tolerance $\epsilon$,
the boundary is split into $n_\Lambda$ disjoint
panels $\Lambda_i$, $i = 1,\ldots,n_\Lambda$ using the adaptive algorithm to be described in Section~\ref{sc:amr}.
Each panel will have $p$ nodes,
giving $N =pn_\Lambda$ nodes in total, hence $2N$ density unknowns.
The $i$th panel is $\Lambda_i = Z([a_{i-1},a_i])$, where
$a_i$, $i=0,\dots,n_\Lambda$
are the parameter breakpoints of all panels
(where $a_{n_\Lambda} \equiv a_0 (\mbox{mod } 2\pi)$).
Let the $p$-point Gauss--Legendre nodes and weights on the parameter
interval $[a_{i-1},a_i]$ be $t^{(i)}_j$ and $w^{(i)}_j$ respectively,
for $j=1,\dots,p$.
Then, for smooth functions $f$ on $\Gamma$, the quadrature rule
\be
\int_\Gamma f(y) ds_y =
\int_0^{2\pi} f(Z(t)) \,|Z'(t)|\, dt
\; \approx \;
\sum_{i=1}^{n_\Lambda} \sum_{j=1}^p f(Z(t^{(i)}_j))\, |Z'(t^{(i)}_j)|\, w^{(i)}_j
=:
\sum_{\ell=1}^N f(y_\ell)\, |Z'(t_\ell)|\, w_\ell
\label{quadr}
\ee
holds to high-order accuracy.
In the last formula above $\{t_\ell\}_{\ell=1}^N$
denotes the entire set of parameter nodes,
$y_\ell=Z(t_\ell)$ their images on $\Gamma$,
and $w_\ell$ their corresponding weights.

The Nystr\"om method \cite[Sec.~12.2]{kress-bie-book} is then used to
approximately solve the BIE \eqref{biesimple}.
Broadly speaking, this involves
substituting the quadrature rule \eqref{quadr} for the integral
implicit in the BIE,
then enforcing the equation at the quadrature nodes themselves.
The result is the $2N$-by-$2N$ linear system,
\be
\bm A \bm\sigma = \mathbf{f}
\label{linsys}
\ee
where $\mathbf{f} := \{f(t_\ell)\}_{\ell=1}^N$
is the vector of 2-component values of the right-hand side
$f$ at the $N$ nodes, and
$\bm\sigma := \{\sigma_\ell\}_{\ell=1}^N$ is the vector of 2-component densities
at the $N$ nodes.
$\bm A = \{A_{ij}\}_{i,j=1,\dots,n_\Lambda}$ is a $n_\Lambda\times n_\Lambda$ block matrix, where each block $A_{ij}$ is a $2p\times2p$ submatrix that represents the interaction from source panel $\Lambda_j$ to target panel $\Lambda_i$.
For targets that are ``far''  (in a sense discussed below)
from a given source panel, the formula for filling corresponding elements
of $\bm A$ is simple.
Letting the index be $\ell$ for such a target node, and $\ell'$ for a source
lying in such a panel,
using the smooth rule \eqref{quadr} 
gives the matrix element
\be
   {\bm A}_{\ell,\ell'} \; = \; K(y_\ell,y_{\ell'})\, |Z'(t_{\ell'})|\, w_{\ell'}
   ~,
   \qquad\mbox{(Nystr\"om rule for $y_\ell$ ``far'' from panel containing $y_{\ell'}$)}
\label{plain}
\ee
where $K(x,y) = S(x,y)+D(x,y)$, the latter being the kernels
appearing in \eqref{sto_slp}--\eqref{sto_dlp}.
Recall that each element in \eqref{plain} is a $2\times 2$ tensor, since
$K$ is.
This defines the plain (smooth) rule for matrix elements (note that we need
not include the diagonal $I/2$ from \eqref{biesimple} since
$\ell=\ell'$ is never a ``far'' interaction).

Sections \ref{s:close_corr}--\ref{s:matel} will be devoted to
defining ``close'' vs ``far'' and
explaining how ``close'' matrix elements are filled.
Assuming for now that this has been done,
the result is a dense matrix $A$ that is well-conditioned independent of $N$,
because the underlying integral equation is of Fredholm second kind.
Then an iterative solver for \eqref{linsys}, such as GMRES,
will converge rapidly.
The result is the vector $\bm\sigma$ approximating the density at the
nodes.
Assuming that \eqref{linsys} has been solved exactly, there
is still a discretization error,
whose convergence rate is known to inherit that of the quadrature scheme
applied to the kernel \cite[Sec.~12.2]{kress-bie-book}.
Given this, the density function may be evaluated at any point
$y\in\Gamma$ using
either the Nystr\"om interpolant 
(which is global and hence inconvenient),
or local $p$th-order Lagrange interpolation
from just the points on the panel in which $y$ lies.
When needed, we will use the latter.
The generalization of the above Nystr\"om method to multiple
closed curves is clear.

\begin{remark}
  When the problem size is large, the matrix-vector multiplication in \eqref{linsys}, which requires $\bigO(N^2)$ time, can be rapidly computed using the fast multipole methods (FMM) in only $\bigO(N)$ time.
  This is because all but $\bigO(N)$ of the matrix elements involve the
  plain rule \eqref{plain}, for which applying the matrix is equivalent
  to evaluation of a potential with weighted source strengths.
  In our large examples (Figures \ref{fig:microfluidic_chip} and \ref{fig:vasnet}), we use an OpenMP-parallelized Stokes FMM code due to Rachh, built upon the Goursat representation of the biharmonic kernel \cite{bhfmm2d,greenbaum92}.
\label{r:fmm}
\end{remark}

Finally, once $\bm \sigma$ has been solved for, the evaluation of
the velocity $u(x)$ at arbitrary targets $x\in\Omega$ is possible,
by approximating the representation
\eqref{urep} or \eqref{urepext}, as appropriate.
By linearity, this
breaks into a sum of contributions from each source panel on each
curve, which may then be handled separately.
Thus a given target $x\in\Omega$ may, again,
fall ``far'' or ``close'' to a given
source panel (denoted by $\Lambda = Z([a,b])$).
If it is ``far'' (according to the same criterion as for matrix elements),
then a simple plain rule is used.
Letting $u_\Lambda$ be the contribution to $u$ from source panel $\Lambda$,
this rule arises, as with \eqref{plain}, simply by substituting
\eqref{quadr} into the representation, giving
\be
u_\Lambda(x) =
\int_\Lambda K(x,y) \sigma(y) ds_y
\approx \sum_{j=1}^p K(x,y_j) \,|Z'(t_j)|\, w_j
\sigma_j~,
\qquad\mbox{(evaluation rule for $x$ ``far'' from $\Lambda$)}
\label{plaineval}
\ee
where $y_j:=Z(t_j)$ are the nodes,
and $\sigma_j:=\sigma(y_j)$ the 2-component
density values, belonging to $\Lambda$.
For large $N$, the FMM is ideal for the task of evaluating
$u$ at many targets, using this plain rule.
An identical quadrature rule may be applied
to the representations in the Appendix to evaluate
pressure and traction at $x$.

\subsection{Close-evaluation and self-evaluation corrections}
\label{s:close_corr}

If a target (on-surface node or off-surface evaluation point) is close enough to a source panel that \eqref{plain} is inaccurate,
a new \emph{close-evaluation} formula is needed.
A special case is when the target is a node from the source panel
itself, which we call \emph{self-evaluation}.
The integrand is nearly singular for close-evaluation and singular for self-evaluation; in the latter case, the integral is understood as an improper integral (SLP) or a principal value integral (DLP).
Conveniently, we will use the same rule for self-evaluation
as both on-surface and off-surface close-evaluation
(this contrasts much prior work, where the self-evaluation used
another set of tailored high-order integration schemes; see, e.g., \cite{hao}).

We quantify ``close'' and ``far'' as follows. Given a panel $\Lambda_j = Z([a,b])$, a target point $x$ is \emph{close} to $\Lambda_j$ if it lies inside the ellipse
\begin{equation}
|z-Z(a)| + |z-Z(b)| = C \, S, \label{eq:panel_close_targ}
\end{equation}
otherwise $x$ is \emph{far} from $\Lambda_j$.
A panel $\Lambda_i$ is close to $\Lambda_j$ if any point $z\in\Lambda_i$ lies inside the ellipse \eqref{eq:panel_close_targ}, otherwise $\Lambda_i$ is far from $\Lambda_j$. (Note that $\Lambda_j$ is close to itself.)
In \eqref{eq:panel_close_targ}, $S$ is the arc length of $\Lambda_j$ and $C > 1$ is a constant.
For the numerical examples in Section \ref{sc:results} we have picked $C=2.5$, which is large enough to include all of both neighboring panels of $\Lambda_j$ most of the time.

The rationale for the above heuristic is based upon the accuracy of the
plain rule \eqref{plaineval} (and its corresponding matrix element rule
\eqref{plain}).
Examining \eqref{quadr}, if
the integrand $f(Z(t))|Z'(t)|$ is analytic with respect to $t$
within a Bernstein ellipse (for the parameter domain $[a_{j-1},a_j]$
for this panel) of size parameter $\varrho>1$,
then the error convergence for the
Gauss--Legendre rule for this panel is $\bigO(\varrho^{-2p})$,
i.e.\ exponential in the panel degree $p$ \cite[Theorem~19.3]{trefethen2013approximation}.
Since in our case $f(Z(t)) = K(x,Z(t))\sigma(Z(t))$, and $K(x,y)$ is analytic
for $x\neq y$, for such analyticity of the integrand,
$x$ must be outside the {\em image} of this Bernstein ellipse under $Z$.
In the case where the panel is approximately flat,
this image is approximated by the ellipse with foci $Z(a)$ and $Z(b)$,
which gives the above geometric ``far'' criterion.
The choice of $C$ is empirically made to achieve an exponential error
convergence rate in $p$ no worse than that due to the next-neighboring
panels discussed in Stage 1 of Section~\ref{sc:amr},
in the case of panel shapes produced by the procedure in that section.
Note that $\sigma(Z(t))$ must also be assumed to be analytic in the
ellipse; we expect this to hold again because the panels will be
sufficiently refined.
For a more detailed error analysis of the plain panel rule
using the Bernstein ellipse, see \cite[Sect.~3.1]{AQBX}.

So far we have presented (only in the ``far'' case) formulae
for both filling Nystr\"om matrix elements \eqref{plain},
and for evaluation of the resulting velocity potential \eqref{plaineval}.
We now make the point that, in both the ``far'' and ``close'' cases,
these are essentially the same task.

\begin{remark}[Matrix-filling is potential evaluation]
The matrix element formula \eqref{plain} is just a special case of
the evaluation formula \eqref{plaineval} for the on-surface target
$x=y_\ell$, and with a Kronecker-delta density $\sigma_j = \delta_{j,\ell'}$.
I.e., one can compute (the ``far'' contributions to) $\bm A\bm\sigma$
from $\bm\sigma$, as needed for each iteration in the solution of
\eqref{linsys}, simply via the evaluation formula \eqref{plaineval}
with targets $x$ as the set of nodes $\{y_\ell\}$.
This will also apply to the special ``close'' formulae presented below.
Henceforth we now present only formulae for evaluation;
the corresponding matrix element formulae are easy to extract
(see Section~\ref{s:matel})
\end{remark}

Note that the diagonal blocks $A_{ii}$ when computed using the special close-evaluation quadratures will implicitly include the $I/2$ jump relation appearing in \eqref{biesimple}.

Finally, to accelerate the computation, the close- and self-evaluations can be precomputed as matrices
(see Section~\ref{s:matel})
from which the (inaccurate) matrices involving the plain rule \eqref{plain}
are subtracted.
The resulting ``correction matrix'' blocks
are assembled and stored as a $2N$-by-$2N$ sparse matrix.
The entire application of $\bm A$ to the density vector
is then performed using the FMM with the plain rule \eqref{plain},
plus the action of this sparse matrix to replace the ``close'' interactions
with their accurate values.
This application is used to solve the whole linear system iteratively.
We do this for our large-scale examples in Section~\ref{sc:results}.

\subsection{Close-evaluation of potentials}
\label{s:eval}

Since we will perform all Nystr\"om matrix filling using the
same formulae as for close-evaluation of potentials (on- or off-surface),
we now present formulae for the close-evaluation task.
As before, we consider a single target point $x\in\Omega$
which is ``close'' to the single source panel $\Lambda = Z([a,b])$,
on which a density $\sigma_j$ is known at the nodes $j=1,\dots,p$.
Recall that the $p$-point Gauss--Legendre nodes and weights for the
parameter interval $[a,b]$ are $t_j$ and $w_j$.

We adapt special panel quadratures proposed for the Laplace case
by Helsing and Ojala \cite{helsing-ojala}.
They use high-order polynomial interpolation in the complex plane to
approximate the density function, then apply a recursion to exactly evaluate
the near-singular integral for each monomial basis function.
In Section~\ref{sc:problem}
we showed that all the needed Stokes potentials may be written
in terms of $I_L(x)$, $I_C(x)$, $I_H(x)$ and $I_S(x)$ from
\eqref{IL}--\eqref{IS}, involving
scalar Cauchy densities $\tau$ derived from the given Stokes density $\sigma$.
Thus in the following subsections we need only cover close-evaluation of
each contour integral in turn.
Although much of this is known \cite{helsing-ojala},
there are certain implementation details and choices
that make a complete distillation valuable.

It turns out that the
Helsing--Ojala polynomial approximation is most stable
when assuming that $Z(a) = 1$ and $Z(b)=-1$, i.e. the panel endpoints
are $\pm 1$ in the complex plane.
Thus we start by making this assumption, then in Sec.~\ref{s:scale}
show how to correct the results for a panel with general endpoints.

Recall that $\sigma$, hence the derived scalar functions $\tau$ needed
in the contour integrands, is available only at the $p$ nodes of $\Lambda$.
In order to improve the accuracy of the complex approximation
for bent panels,
firstly an upsampling is performed to $m>p$ ``fine'' nodes, using 
Lagrange interpolation in the parameter $t\in [a,b]$ from the $p$ nodes
to the $m$ fine nodes.
We find that $m=2p$ is beneficial without incurring significant extra cost.
Let $\tilde\tau_j$, $j=1,\dots,m$ denote the fine density values,
and $z_j=Z(\tilde t_j) \in\Lambda$ be the fine nodes,
where $\tilde t_j$ and $\tilde w_j$ are the $m$-point Gauss--Legendre nodes and
weights respectively for $[a,b]$.
The following schemes now will use only the fine values and nodes.

\subsubsection{Close evaluation of the Cauchy potential}
\label{sec:IC}

We approximate $\tau$ on the panel $\Lambda$ in the complex variable by
the degree $m-1$ polynomial
\begin{equation}
\tau(y) \approx \sum_{k=1}^m a_k y^{k-1},\qquad y\in\Lambda ~.
\label{poly}
\end{equation}
The vector of coefficients $\mathbf{a} := \{a_k\}_{k=1}^m$
is most conveniently found by using a dense
direct solve of the square Vandermonde system
\be
V \mathbf{a} = \tilde{\bm\tau}~,
\label{vand}
\ee
with elements $V_{jk} = z_j^{k-1}$, $j,k=1,\dots,m$, and right hand side
$\tilde{\bm{\tau}} := \{\tilde\tau_j\}_{j=1}^m$.
It is known that, for
any arrangement of nodes $z_j$ other than those very close to the
roots of unity,
the condition number of $V$ grows exponentially with $m$ \cite{pan16}.
However, as discussed in \cite[App.~A]{helsing-ojala},
at least for $m<50$, despite the extreme ill-conditioning,
{\em backward stability} of the linear solver insures that
the resulting polynomial matches the values at the nodes to
close to machine precision.
For this we use MATLAB {\tt mldivide}
which employs standard partially-pivoted Gaussian elimination.

\begin{figure}[t]
\centering
\begin{tabular}{cc}
(a) & (b)\\
\includegraphics[width=0.25\textwidth]{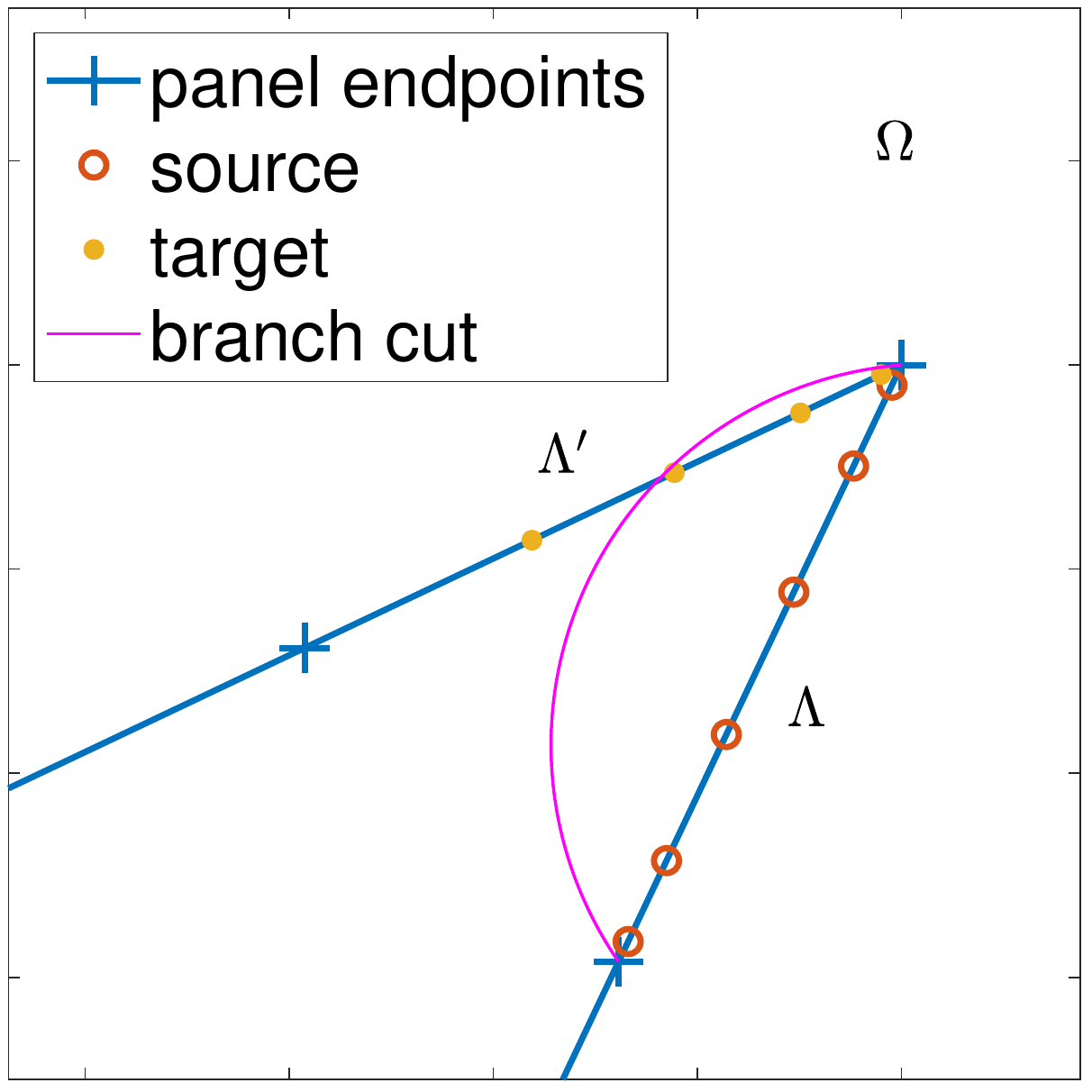} & \includegraphics[width=0.25\textwidth]{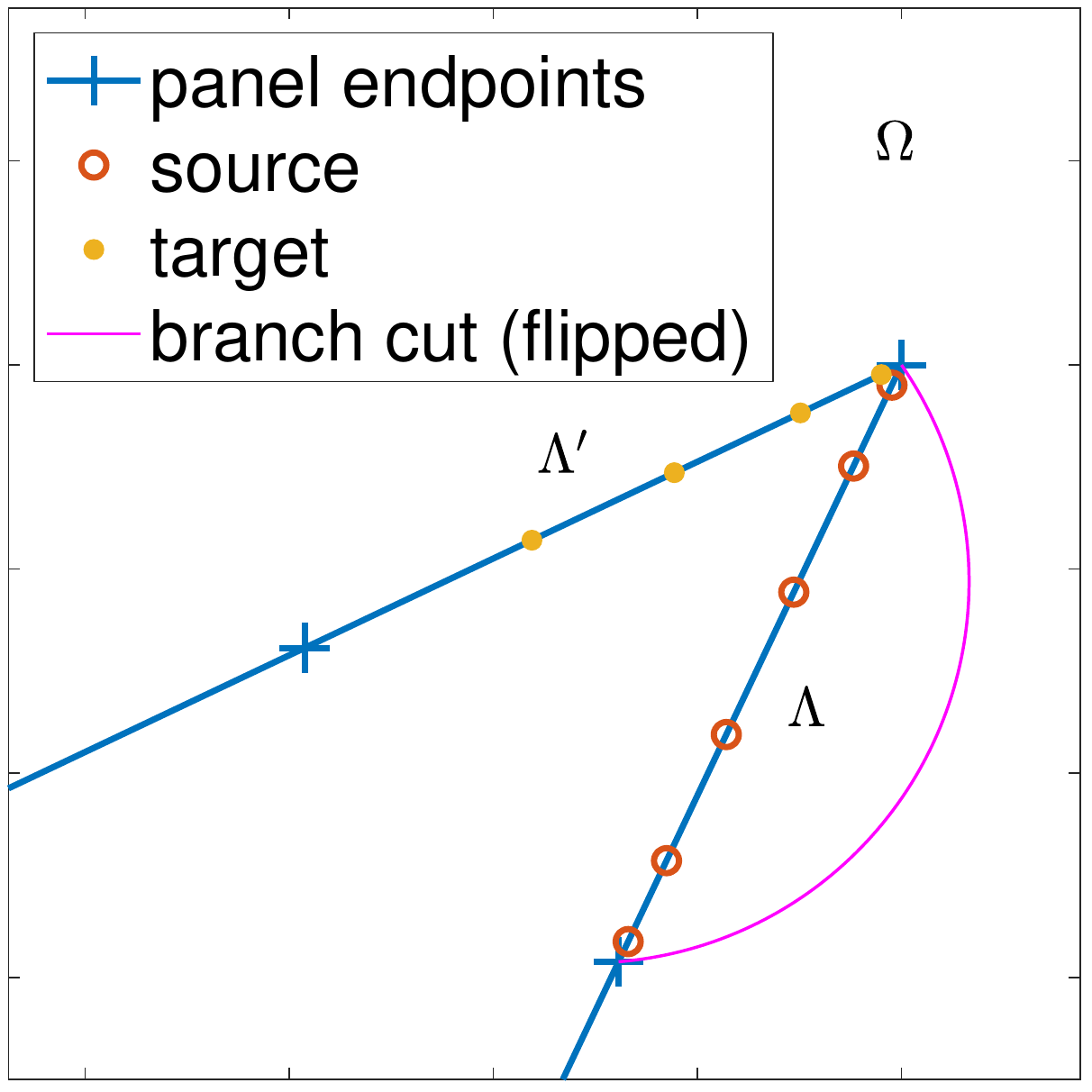}
\end{tabular}
\caption{Special handling of close evaluation branch cut when the panel is touching a reentrant corner. (a) The target panel $\Lambda'$ is crossing the branch cut of the source panel $\Lambda$ defined by \eqref{p1use}, resulting in wrong close evaluation values. (b) Changing the sign of $\phi$ in \eqref{p1use} flips the branch cut to the other side; close evaluation at the targets are now correct.}\label{fig:corner_branchcut}
\end{figure}

The remaining step is to compute the contour integral of each
monomial,
\begin{equation}
p_k = p_k(x) := \int_{-1}^1 \frac{y^{k-1}}{y-x}dy,\quad k=1,\ldots,m~,
\label{eq:p_k}
\end{equation}
which are, recalling \eqref{IC}, then combined using \eqref{poly} to get
\begin{equation}
I_C \;\approx\; \sum_{k=1}^m a_k p_k~.
\label{ICsum}
\end{equation}
By design, since the monomials are with respect to $y$ in the complex
plane (as opposed to, say, the parameter $t$),
Cauchy's theorem implies that each $p_k$ is independent of the curve $\Lambda$
and depends only on the end-points.
Specifically, for $k=1$ we may integrate analytically by deforming $\Lambda$
to the real interval $[-1,1]$, so
\be
p_1 \;:=\; \int_{-1}^1 \frac{1}{y-x}dy
\;=\;
\log(1-x)-\log(-1-x)\pm 2\pi i {\mathcal N}_x
\label{p1}
\ee
where ${\mathcal N}_x \in \mathbb{Z}$ is an integer winding number that
depends on the choice of branch cut of the log function.
For the standard cut on the negative real axis then
${\mathcal N}_x=0$ when $x$ is outside the domain enclosed by
the oriented curve composed of $\Lambda$ traversed forwards plus $[-1,1]$
traversed backwards, ${\mathcal N}_x=+1$ ($-1$) when $x$ is inside a region
enclosed counterclockwise (clockwise) \cite{helsing-ojala}.
However, since it is inconvenient and error-prone to decide ${\mathcal N}_x$
for points on or very close to $\Lambda$ and $[-1,1]$, we prefer
to combine the two logs then effectively
rotate its branch cut
by a phase $\phi\in\mathbb{R}$, by using
\be
p_1 \;=\; i\phi + \log \frac{1-x}{e^{i\phi}(-1-x)}~,
\label{p1use}
\ee
where $\phi = -\pi/4$ when the upwards normal of the panel points into
$\Omega$ (as for an interior curve),
or $\phi=\pi/4$ otherwise.
This has the effect of pushing the branch cut ``behind'' the panel
(away from $\Omega$; see Fig.~\ref{fig:corner_branchcut}),
with the cut meeting $\pm 1$ at an angle $\phi$
from the real axis.
The potential is correct in the closure of $\Omega$, including
on the panel itself, without any topological tests needed (hence the unified handling of close and self evaluations in Sec.~\ref{s:close_corr}).
This can fail if a panel is very curved (such a panel would be
inaccurate anyway), or if a piece of $\Omega$
approaches close to the back side of the panel (which can be prevented
by adaptive refinement as in Section~\ref{sc:amr}).
In practice we find that this is robust. However, we will mention one special situation where \eqref{p1use} could fail and therefore careful adjustment of the branch cut is critical; see Remark \ref{rmk:corner_branchcut} below.

The following 2-term recurrence is easy to check by adding and subtracting
$xy^{k-1}$ from the numerator of the formula \eqref{eq:p_k} for $p_{k+1}$:
\begin{equation}
p_{k+1} = xp_k + (1-(-1)^k)/k~.
\label{eq:p_k_recur}
\end{equation}
For $|x|<1.1$ we find that the recurrence is sufficiently stable to use
upwards from the value $p_1$ computed by \eqref{p1use}, to get $p_2,\dots,p_m$.
However, for targets outside this close disc, especially for larger $m$,
the upwards direction is unstable.
Thus, here instead we use numerical quadrature on \eqref{eq:p_k} to get
\be
p_m \approx \sum_{j=1}^m \frac{z_j^{m-1}}{z_j - x} Z'(\tilde t_j) \tilde w_j~,
\label{pm}
\ee
then apply \eqref{eq:p_k_recur} downwards to compute $p_{m-1},\dots,p_2$,
and as before use $p_1$ from \eqref{p1use}. Outside of the disc, no
branch cut issues arise.

\begin{remark}[branch cuts at corners]\label{rmk:corner_branchcut}
When a panel $\Lambda$ is directly touching a corner, directly applying \eqref{p1use} can fail no matter how much the panels are refined. In Figure \ref{fig:corner_branchcut}a, the panel on the opposite side of the corner, $\Lambda'$, is always behind $\Lambda$ and lying across the branch cut associated to $\Lambda$. Consequently, the close evaluation from $\Lambda$ to $\Lambda'$ results in completely wrong values, also leading to ill-conditioning of the whole system \eqref{linsys}. Instead, one can simply change the sign of $\phi$ in \eqref{p1use} to flip the branch cut to accommodate the targets on $\Lambda'$ (Figure \ref{fig:corner_branchcut}b). In practice, we find that this is robust for corners of arbitrary angles.
\end{remark}

\subsubsection{Close evaluation of the logarithmic potential}
\label{sec:IL}

For $\tau$ real, we can write the final form in \eqref{IL} as
\be
I_L \;=\; \re \int_\Lambda \log (y-x) \frac{\tau(y)}{in_y} dy,
\label{ILre}
\ee
which shows that the quantity to approximate on $\Lambda$
as a complex polynomial is $\tau(y)/in_y$. Thus we find the coefficients in
\be
\frac{\tau(y)}{in_y} \;\approx\; \sum_{k=1}^m b_k y^{k-1}, \qquad y\in\Lambda~,
\label{polyb}
\ee
by solving \eqref{vand} but with modified right hand side
$\{\tilde \tau_j/i n_{z_j}\}_{j=1}^m$.
Defining
\begin{equation}
q_k \;:=\; \int_{-1}^1 y^{k-1}\log (y-x) dy,\quad k=1,\ldots,m
\end{equation}
then \eqref{polyb} gives
\begin{equation}
\int_{-1}^{1} \log (y-x) \frac{\tau(y)}{in_y}dy \;\approx\;
\sum_{k=1}^m b_k q_k~,
\label{eq:tau_by_n}
\end{equation}
whose real part is $I_L$.
Each $q_k$ is computed from $p_k$ of \eqref{eq:p_k}
as evaluated in Sec.~\ref{sec:IC}, via a formula easily proven
by integration by parts,
\begin{equation}
q_{k} = \frac{-p_{k+1} + \log(1-x) -(-1)^k \log(1+x)}{k}
=\begin{cases}
\bigl(-p_{k+1} + i\phi + \log\frac{1-x}{e^{i\phi}(-1-x)} \bigr) / k~,&
k \mbox{ even },\\
\bigl(-p_{k+1} + \log[(1-x)(-1-x)] \bigr) / k~,&
k \mbox{ odd },\\
\end{cases}
\label{qk}
\end{equation}
where the latter form is that used in the code, needed to match the
branch cut rotation used for $p_k$ in \eqref{p1use}.
The final evaluation of $I_L$ is then via
\be
I_L \;\approx\; \re \sum_{k=1}^m b_k q_k~.
\ee

\subsubsection{Close evaluation of the Hadamard and supersingular potentials}
\label{sec:IH}

The double-layer Stokes velocity requires 
gradients of Laplace potentials
\eqref{gradD}, which require $I_H$.
Also, the traction of the Stokes single-layer \eqref{tracIH} involves $I_H$
applied to $\tau(y)/in_y$, and the traction of the Stokes double-layer
further involves $I_S$ (Appendix~\ref{a:trac}).

Using the complex monomial expansion \eqref{poly}, we have
\begin{equation}
  I_H \;\approx \; \sum_{k=1}^m a_k r_k, \hspace{1in}
  I_S \;\approx \; \sum_{k=1}^m a_k s_k,
\end{equation}
where
\begin{equation}
  r_k \;:=\; \int_{-1}^1 \frac{y^{k-1}}{(y-x)^2}dy, \hspace{.5in}
    s_k \;:=\; \int_{-1}^1 \frac{y^{k-1}}{(y-x)^3}dy,
  \qquad k=1,\ldots,m~.
\label{eq:r_k}
\end{equation}
The following formulae can be shown by integration by parts, and
enable $r_k$ and $s_k$ to be found,
\begin{equation}
  r_{k} = (k-1)p_{k-1} + \frac{(-1)^{k-1}}{-1-x} - \frac{1}{1-x}~,
  \qquad
  s_{k} = \frac{k-1}{2}r_{k-1} + \frac{(-1)^{k-1}}{2(-1-x)^2} - \frac{1}{2(1-x)^2}~,
\label{eq:r_k_recur}
\end{equation}
using $p_k$ from \eqref{eq:p_k} as computed in Sec.~\ref{sec:IC}, and $p_0=0$.

\subsubsection{Transforming for general panel endpoints}
\label{s:scale}

To apply close-evaluation methods in the above three sections to a general panel
$\Lambda = Z([a,b])$,
define the complex scale factor $s_0 := (Z(b)-Z(a))/2$ and
origin $x_0 = (Z(b)+Z(a))/2$.
Then the affine map
$$
x = s(\tilde x) := \frac{\tilde x- x_0}{s_0}
$$
takes any target $\tilde x$ to its scaled version $x$.
Likewise, the fine nodes are scaled by $z_j = s(Z(\tilde t_j))$,
and the factor $Z'$ in \eqref{pm} is replaced by $Z'/s_0$.
Following Sec.~\ref{sec:IC} using these scaled target and fine nodes,
no change in the result $I_C$ is needed.
However, the value of $I_H$ computed in
Sec.~\ref{sec:IH} must afterwards be divided by $s_0$,
and the value of $I_S$ divided by $s_0^2$.
The value of $I_L$ computed in
Sec.~\ref{sec:IL} must be multiplied by $|s_0|$,
and then have $|Z'(\tilde t_j) w_j/s_0| \log |s_0|$
subtracted.

\subsection{Computation of close-evaluation matrix blocks}
\label{s:matel}

The above described how to evaluate
$(I_L\tau)(x)$,
$(I_C\tau)(x)$,
$(I_H\tau)(x)$ and
$(I_S\tau)(x)$ given known samples $\tilde\tau_j$ at a panel's fine nodes.
In practice it is useful to instead precompute a matrix block $A$ which
takes any density values at a panel's original $p$ nodes $y_j$ to
a set of $n$ target values of the contour integral.
Consider the case of the Cauchy kernel,
and let $A$ denote this $n$-by-$p$ matrix.
Let $L$ be the $m$-by-$p$ Lagrange interpolation matrix from the
nodes $t_j$ to fine nodes $\tilde t_j$, which need be filled once and for all.
Let $P$ be the $n$-by-$m$ matrix with entries $P_{ik} = p_k(x_i)$,
given by \eqref{eq:p_k},
where $\{x_i\}_{i=1}^n$ is the set of desired targets.
In exact arithmetic one has
$$
A = P V^{-1} L~.
$$
However, since $V$ is very ill-conditioned,
filling $V^{-1}$ and using it to multply to the right is numerically unstable.
Instead an {\em adjoint method} is used:
one first solves the matrix equation $V^\top X = P^\top$, where $\top$ indicates
non-conjugate transpose, then forms the product
$$
A =  X^\top L~.
$$
The matrix solve is done in a backward stable fashion via MATLAB's
{\tt mldivide}.
A further advantage of the adjoint approach is that if $n$ is small,
the solve is faster than computing $V^{-1}$.

The formulae for the logarithmic, Hadamard and supersingular kernels are analogous.


\section{Adaptive panel refinement}\label{sc:amr}

In order to solve a BIE to high accuracy, it is necessary to set up panels such that the given complex geometry is correctly resolved. In this section, we describe a procedure that adaptively refine the panels based solely on the geometric properties. Specifically, our refinement algorithms take into account the accuracy of geometric representations (including arc length and curvature), the location of corners, and the distance between boundary components.
It necessarily has some ad-hoc aspects, yet we find it quite robust in practice.

Suppose that for a user-prescribed tolerance $\epsilon$, the goal is to find a partition $\Gamma = \bigcup\limits_{i=1}^{n_\Lambda}\Lambda_i$ such that the error, $\varepsilon$, of evaluating boundary integral operators such as \eqref{biesimple} satisfies $\varepsilon \lessapprox \epsilon$. To this end, we describe our adaptive refinement scheme which proceeds with three stages. In what follows, we again assume that the panel under consideration $\Lambda=Z([a,b])$ is rescaled such that its two endpoints are $\pm1$.

\paragraph{Stage 1: Choice of overall $p$.} Given tolerance $\epsilon$, the goal is to determine a number of quadrature points, $p$, applied to all panels, such that the relative quadrature error on any panel is $\bigO(\epsilon)$.
As mentioned above, the $p$-point Gauss--Legendre quadrature on $[-1,1]$ has $\bigO(\varrho^{-2p})$ error if the integrand can be analytically extended to a Bernstein ellipse of parameter $\varrho>1$, where the semi-major axis of this ellipse is $(\varrho+\varrho^{-1})/2$ \cite[Thm.~19.3]{trefethen2013approximation}. Therefore, making $\varrho^{-2p} \leq \epsilon$ we obtain the first term of the right-hand side in 
\begin{equation}
p \geq \left\lceil\frac{\log_{10}(1/\epsilon)}{2\log_{10} \varrho}\right\rceil
+ c, \label{eq:pan_order}
\end{equation}
where the second term accounts for unknown prefactors.
Empirically we set $c=1$.

To determine $\varrho$, we require that the Bernstein $\varrho$-ellipse of each panel encloses both its immediate neighboring panels.
This insures that, when applying the smooth quadrature rule \eqref{quadr} between the nearest non-neighboring (``far'') panels that do not touch a corner,
the integrand continues to a function analytic inside the
$\varrho$-ellipse, so, by the above discussion, the relative error is no worse than $\epsilon$.
(This will not apply to panels touching a corner,
but they are small enough to have negligible contributions.)
Stages 2--3 below will place an upper bound of $\lambda$ on the ratio of the
lengths (with respect to parameter) of neighboring panels.
Combining these two relations gives
\begin{equation}
\frac{\varrho+\varrho^{-1}}{2} = 1 + \frac{2}{\lambda}.\label{eq:find_rho}
\end{equation}
One then solves \eqref{eq:find_rho} for $\varrho$ and substitutes it into \eqref{eq:pan_order} to obtain a lower bound for $p$. For example, $\lambda\le 3$ holds in our examples, so $\varrho = 3$, and therefore
we have as sufficient the simple rule $p = \lceil\log_{10}(1/\epsilon)+1\rceil$.

\paragraph{Stage 2: Local geometric refinement.} In this stage, panels are split based on local geometric properties:

\begin{enumerate}
\item {\em Corner refinement}. Panels near a corner are refined
geometrically so that each panel is a factor $\lambda$ shorter in parameter
than its neighbor (see lines 8--11 of Algorithm \ref{alg:loc_refine}).

To each corner is associated a factor $\lambda\geq2$.
A rule of thumb is to use $\lambda = 2$ for sharper corners (e.g.\ whose angle $\theta$ is close to $0$ or $2\pi$) which are harder to resolve,
and use $\lambda > 2$ for ``flatter'' corners (e.g.\ $\theta$ closer to $\pi$) to reduce the number of panels without affecting the overall achieved accuracy.
In practice, we use $\lambda = 3$ for corners $\pi/2\leq\theta\leq3\pi/2$; 
for a problem with many flat corners, this can reduce the total number of unknowns by a factor of about 2/3 (or about $\log_3 2$).

Near a corner, refinement stops when the panels touching the corner are shorter than $\epsilon^{\alpha}$, where $\alpha$ is an empirical power parameter chosen for each corner.
Recent theoretical results for the plain double-layer formulation for
the Stokes Dirichlet BVP state that the density is
a constant plus a bounded singular function whose power exceeds $1/2$
for any corner angle in $(0,2\pi)$
\cite{rachh2017solution};
this is similar to the Laplace case \cite{zargaryan84}.
For our $\mathbf{D}+\mathbf{S}$ formulation we observe a density behavior consistent with this.
This might suggest choosing $\alpha=1$ for any corner angle.
In fact, for small (non-reentrant) angles we are able to reduce $\alpha$
somewhat without loss of accuracy, hence do so, to reduce the number of panels.

\item {\em Bent panel refinement}. Panels away from any corners are refined based on how well the smooth geometric properties are represented by the interpolants on their $p$ Legendre nodes. 
We measure the accuracy of geometric representations by the interpolation errors of a set of test functions on a set of test points. 
First, we define the set of test functions $G=\{g_1,g_2,g_3,\ldots\}$ to be approximated on the panel $\Lambda$. The following list of functions are included in $G$ whenever the necessary derivatives are available:
        \begin{itemize}
        \item $g_1(t) = Z(t)$, which resolves the geometry representation.
        \item $g_2(t) = |Z'(t)|$, which resolves arc length, recalling that arc length is $$S = \int_\Lambda ds = \int_a^b |Z'(t)|\,dt$$
        \item $g_3(t) = \frac{|\mathrm{Im}\,(Z''(t)/Z'(t))|^2}{|Z'(t)|}$, which resolves bending, since bending energy is $$E = \int_\Lambda \kappa^2 ds =  \int_a^b \left|\frac{Z'(t)\times Z''(t)}{|Z'(t)|^3}\right|^2|Z'(t)|\,dt = \int_a^b \frac{\{\mathrm{Im}\,(Z''(t)/Z'(t))\}^2}{|Z'(t)|}\,dt$$
        \end{itemize}
    Next, we define the test points to be the $m$ equally spaced points on $[a,b]$, denoted $\tilde t^\Lambda_j, j=1,\ldots,m$,
    and let $t^\Lambda_j$, $j=1,\ldots,p$ be the Legendre nodes. Then for each
    $i=1,2,3$, the relative error of approximating $g_i$ is
    $$\varepsilon_i = \frac{\norm{\tilde{\mathbf{g}}_i - \mathbf{M}\cdot \mathbf{g}_i}}{\norm{\mathbf{g}_i}},$$
    where $\tilde{\mathbf{g}}_i := (g_i(\tilde t^\Lambda_1),\ldots, g_i(\tilde t^\Lambda_m))$, $\mathbf{g}_i := (g_i(t^\Lambda_1),\ldots, g_i(t^\Lambda_{p}))$, and $\mathbf{M}$ is the $m\times p$ interpolation matrix from the Legendre nodes to the test points.
        The panels are refined until $\max_i \varepsilon_i < \epsilon^\beta$, where $\beta>0$ is another power parameter, with default value $\beta = 1$.
\end{enumerate}

The corner and bent panel refinement rules are applied to all panels recursively. The complete procedure is summarized in Algorithm~\ref{alg:loc_refine}.

\begin{algorithm}[t] 
    \caption{Local geometric refinement} \label{alg:loc_refine}
    \algrenewcommand\algorithmiccomment[2][\footnotesize]{{#1\hfill\(\triangleright\) #2}}

    \Input{ The current panel $\Lambda = Z([a,b])$; tolerance $\epsilon$; corner information $C=\{t^c_j, \lambda_j, \alpha_j\}_{j=1}^k$; test function(s) $G = \{g_1,g_2,g_3,\ldots\}$; $\beta$ is the tolerance exponent for the test functions, default to be $1$.}
    
    \begin{algorithmic}[1]
    \Function{Refine}{$Z([a,b]),\epsilon, C, G,\beta$}
        \State Panel parametric length $L = b-a$ 
        \State Panel arc length $S = \int_a^b |Z'|$
        \If{geometry has corners}
            \State Let $t^c_i$ be the corner closest to $[a,b]$
            \If{$L < \epsilon^{\alpha_i}$ or $S < \epsilon^{\alpha_i}$} 
                	\State\Return{$\{a, b\}$} \Comment{panel length reached lower limit, do not split}
            \EndIf
            \If{$t_i^c$ is close enough to the panel $[a, b]$} 
            		\If{$a < t_i^c <b$} $s =  t_i^c$ \Comment{split right at the corner}
    		\ElsIf{$t_i^c < a$} $s = a + {L}/{\lambda_i}$ \Comment{split towards the corner}
    		\ElsIf{$t_i^c > b$} $s = b - {L}/{\lambda_i}$ \Comment{split towards the corner}
    		\EndIf  
            \EndIf
        \EndIf
        \If{split point $s$ is \emph{not} defined}
            \State $\mathbf{g}_i =  g_i$ values at quadrature points
            \State $\tilde{\mathbf{g}}_i = g_i$ values at test points
            \State $\mathbf{M} = $ interpolation matrix from quadrature points to test points
            \State $\varepsilon = \max\limits_{g_i\in G} \norm{\tilde{\mathbf{g}}_i - \mathbf{M}\cdot \mathbf{g}_i}/\norm{\mathbf{g}_i}$ \Comment{interpolation error of test function(s)}
            \If{$\varepsilon>\epsilon^\beta$}
        		  \State $s = {(a+b)}/{2}$ \Comment{split in half}
            \EndIf
        \EndIf
        \If{split point $s$ is defined}
    		  \State $\mathbf{t}^1=$ \Call{Refine}{$Z([a,s]),\epsilon, C, G,\beta$}
    		  \State $\mathbf{t}^2=$ \Call{Refine}{$Z([s,b]),\epsilon, C, G,\beta$}
    		  \State \Return{$\mathbf{t}^1\cup\mathbf{t}^2$}  \Comment{\footnotesize{recursively refine panel}}
	\EndIf
	\State\Return{$\{a, b\}$} \Comment{do not split}
	\EndFunction
        \end{algorithmic}
\end{algorithm}

\paragraph{Stage 3: Global closeness refinement.} At this final stage, panels are further refined if they are (relatively) too close to any non-neighboring panels. Specifically, let $\Lambda^{\text{left}}$ and $\Lambda^{\text{right}}$ be the two immediate neighboring panels of $\Lambda$, and define $\Gamma^{\text{far}}:=\Gamma\setminus(\Lambda^{\text{left}}\cup\Lambda\cup\Lambda^{\text{right}})$ as all non-neighboring panels of $\Lambda$.
Then the panel $\Lambda$ is refined if $d(\Lambda,\Gamma^{\text{far}})$, its distance from $\Gamma^{\text{far}}$, is shorter than its arc length by a factor of 3 (see Line 4-10 of Algorithm~\ref{alg:glob_refine});
see Remark~\ref{r:rootdist_rule} for an alternative, less restrictive, refinement criterion. 

In practice, the distance $d(\Lambda,\Gamma^{\text{far}})$ can be approximated by $\min_{i,j}|y_i - y_j|$, where the minimum is searched among all pairs of nodes $y_i\in\Lambda$ and $y_j\in\Gamma^{\text{far}}$. A kd-tree algorithm \cite{ataiya-kdtree} is used to accelerate this process for our large examples in Section~\ref{sc:results}, in which case the elliptical close neighborhood \eqref{eq:panel_close_targ} is also replaced by $\bigcup_{i=1}^p B(Z_i, C\, S)$, the union of disks around each node on $\Lambda$, for convenience.

The above refinement process is applied to each panel from the output of the previous stage and repeats until no further splitting. The algorithm for this stage is summarized in Algorithm~\ref{alg:glob_refine}. 
We note that since our algorithm is panel-based, it is agnostic of whether two touching panels belong to the same boundary component or not.
Hence this algorithm handles two situations simultanteously: the case of close-touching between different boundary components, as well as the case of ``self-touching'' where a boundary component is almost touching itself.

\begin{algorithm}[t] 
    \caption{Global closeness refinement} \label{alg:glob_refine}
    \algrenewcommand\algorithmiccomment[2][\footnotesize]{{#1\hfill\(\triangleright\) #2}}

    \Input{
	    The refined panels $\Gamma = \bigcup_{k} \Lambda_k$ from Stage 2 (local geometric refinement).
    }
    
    \begin{algorithmic}[1]
     \Function{CloseRefine}{$\Gamma$}
     \State Initialize the output set $\mathbf{t} = \emptyset$, which will contain the final panel endpoints
     \State Initialize the set of new endpoints $\mathbf{t}^\mathrm{new} = \{$endpoints of $\Gamma = \bigcup_{k} \Lambda_k\}$
     \While{$\mathbf{t}^\mathrm{new}\neq\emptyset$}\Comment{repeat until no further splitting}
        \State $\mathbf{t}=\mathbf{t}\cup\mathbf{t}^\mathrm{new}$
        \State $\mathbf{t}^\mathrm{new} = \emptyset$ 
        \State Update panels $\Gamma = \bigcup_{k} \Lambda_k$ based on $\mathbf{t}$ \Comment{ready for a new round of refinement}
     \For{each panel $Z([a,b])\subset\Gamma$}

        \State Locate $\Lambda^{\text{left}}$ and $\Lambda^{\text{right}}$, the two immediate neighboring panels of $Z([a,b])$
        \State Define $\Gamma^{\text{far}}=\Gamma\setminus(\Lambda^{\text{left}}\cup Z([a,b])\cup\Lambda^{\text{right}})$
        \State Compute the distance $d = d(Z([a,b]), \Gamma^{\text{far}})$
        \State Calculate $S = $ arc length of $Z([a,b])$
        \If{$\frac{1}{3}\,S>d$}
            \State $\mathbf{t}^\mathrm{new} = \mathbf{t}^\mathrm{new} \cup \{\frac{a+b}{2}\}$ 	\Comment{split in half}	
        \EndIf
        \EndFor
        \EndWhile
        \State\Return{final panel endpoints $\mathbf{t}$}
	\EndFunction
        \end{algorithmic}
\end{algorithm}

\begin{remark}[Expected convergence rate with corners]  
The above three stages involve two quantities---the panel order $p$,
and a resulting number of panels per corner---both of which grow linearly
with $\log 1/\epsilon$. However, $N$ is the product of these two
quantities, thus, in the presence of corners,
one expects asymptotically $N = \bigO(\log^2 1/\epsilon)$ as
$\epsilon\to0$.
In other words, the error converges {\em root exponentially} in $N$,
i.e.\ as $\bigO(e^{-c \sqrt{N}})$.
This matches the theoretical convergence rate for $hp$-BEM
on polygons by Heuer--Stephan \cite{heuer96}.
This rate has also recently been observed and proven for a geometrically
graded ``method of fundamental solutions'' approach to polygons
by Gopal--Trefethen \cite{rational-corner-laplace}.
\label{r:rootexp}
\end{remark}


\section{Numerical results and discussion}\label{sc:results}
In what follows, numerical examples will be presented to test the overall solution scheme presented so far. In each example, we solve a Dirichlet problem in the domain exterior to the given geometries. The BIE formulation of the problem is $(\frac{1}{2}+\mathbb{S}+\mathbb{D})\sigma = f$ as described in Section \ref{sc:formulation}.

The solution procedure is to first adaptively refine the representation of the geometry using our refinement scheme (Section \ref{sc:amr}), then the BIE is discretized using the special quadrature (Section \ref{sc:panel}) and solved for the density $\sigma$, and finally the solution $u = u_\infty + (\mathbf{S}+\mathbf{D})\sigma$ is evaluated everywhere in the exterior domain on a grid of spacing $\Delta x=0.02$.

We mention that our solution scheme has been tested on boundary value problems
with inhomogeneous boundary data extracted from an analytically known
smooth flow $u$, and, as expected, achieves superalgebraic convergence.
However, in the presence of corners, such smooth test problems do {\em not}
involve the corner singularities that generically arise in physical problems.
For this reason, we only present results on physical flows such as imposed uniform or linear shear flows. In all the examples, the exact solution is not known analytically; therefore, we use the finest grid solution as the reference solution. 

\begin{figure}[t] 
\centering
\begin{tabular}{ccc}
(a) & (b) & (c)\\
\includegraphics[width=0.3\textwidth]{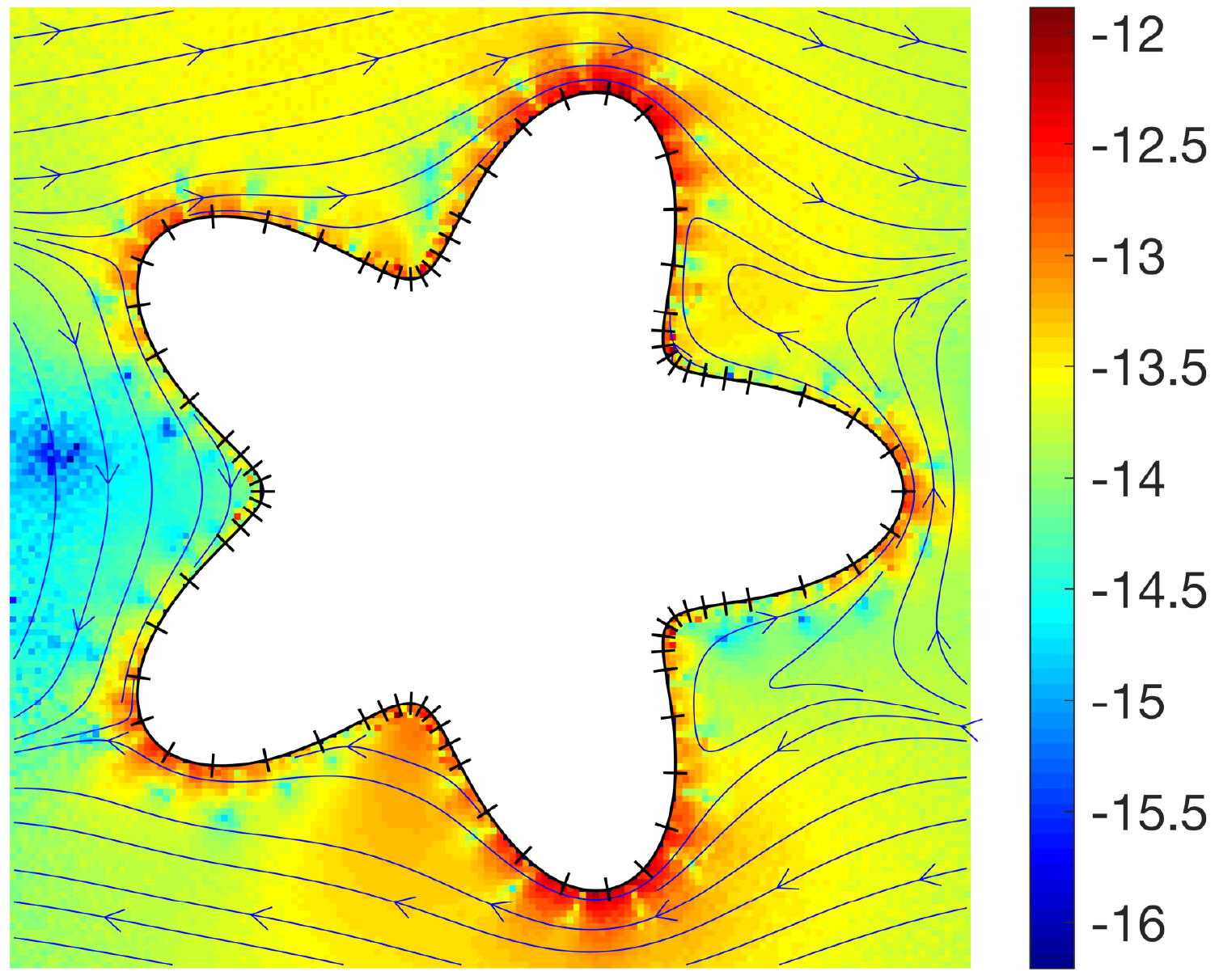}  & \includegraphics[width=0.32\textwidth]{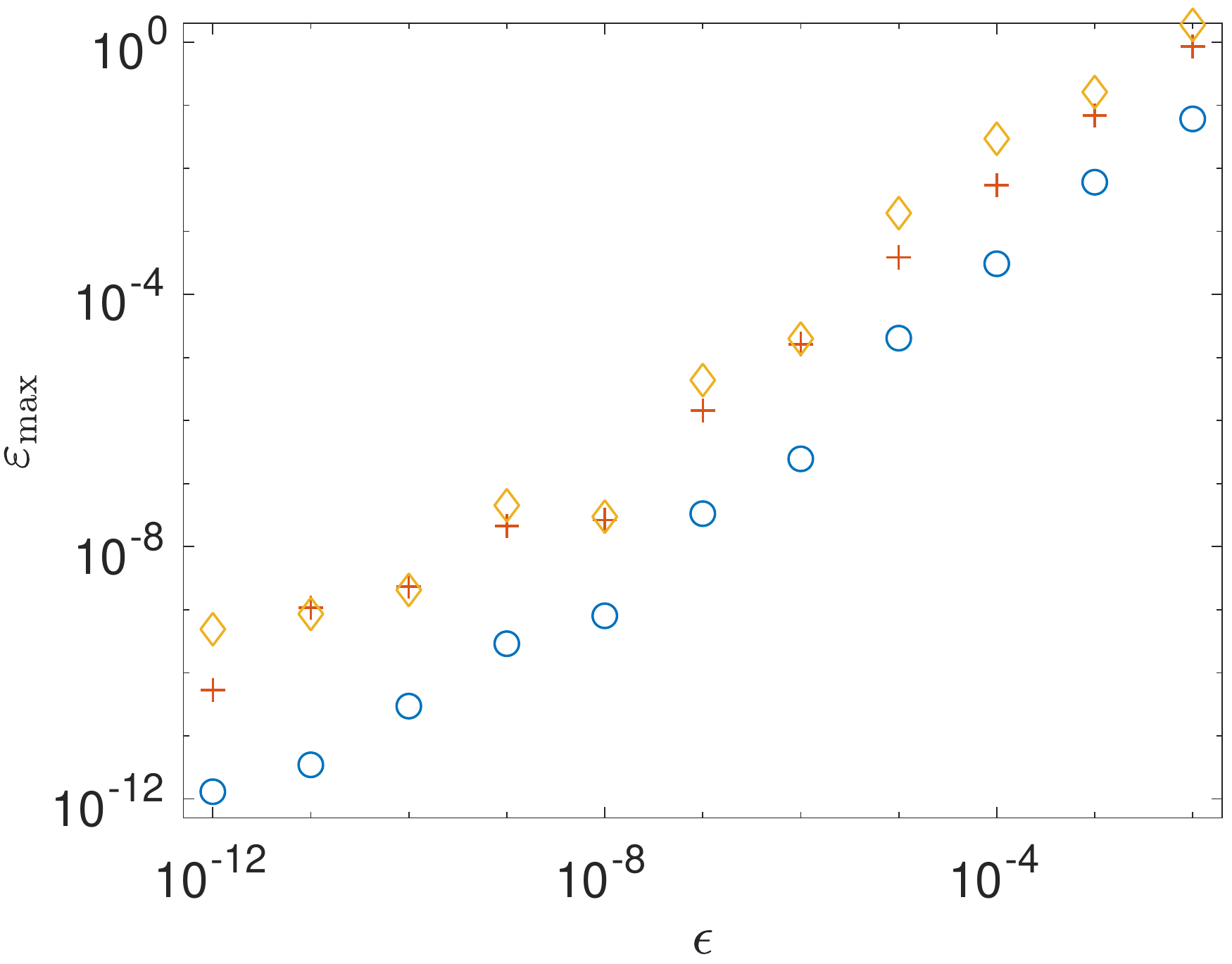} & \includegraphics[width=0.32\textwidth]{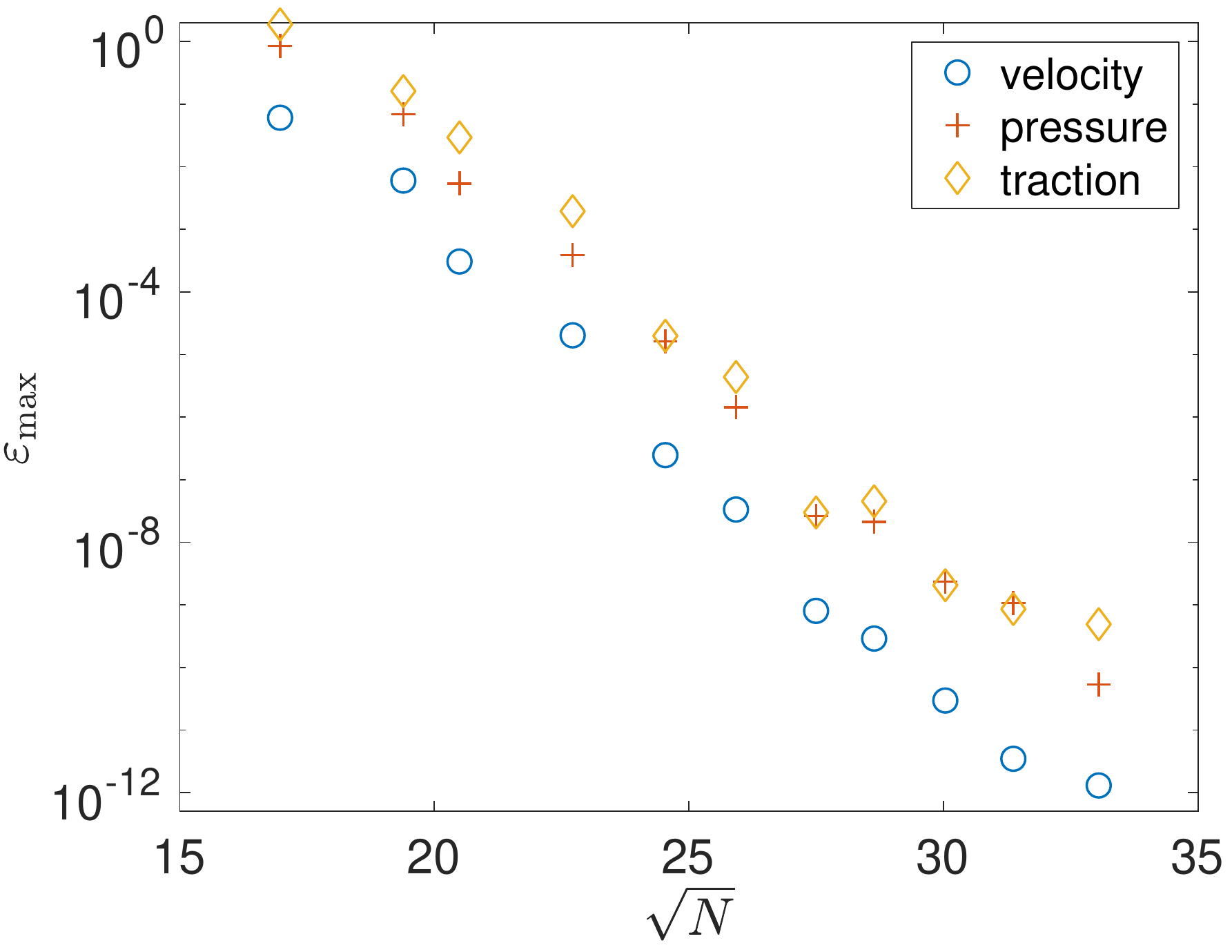}
\end{tabular}
\caption{(a) Linear shear flow past a starfish-shaped island. Streamlines of the flow and panel endpoints (small segments) are shown. Color represents the $\log_{10}$ error of the velocity computed under tolerance $\epsilon=10^{-12}$ (resulted in  $2N=2184$ degrees of freedom).
(b) Convergence of the maximum relative error $\varepsilon_\mathrm{max}$ versus requested tolerance $\epsilon$, where the traction field is computed along the $(1,2)$-direction.
(c) Convergence of the maximum error $\varepsilon_\mathrm{max}$ versus the square root of the number of nodes.
\label{fig:error_ex1}}
\end{figure}

\paragraph{Example 1. Smooth domain.} This example tests our scheme on a linear shear flow around a smooth starfish-shaped island defined by the polar function $r(\theta) = 1 + 0.3 \cos 5\theta$, with no-slip boundary condition $u|_\Gamma=0$ and $u_\infty(x_1,x_2)=(-x_2, 0)$ as $|x|\to\infty$. We have used $\beta = 0.8$ (Line 17 of Algorithm \ref{alg:loc_refine}) for this problem to reduce $N$.
In addition to the velocity field, we have also investigated the convergence of the pressure field and the traction field in the $(1,2)$-direction, both of which are obtained using our close evaluation scheme (Section~\ref{s:eval} and Appendix~\ref{a:trac}). Our scheme achieved accuracies that are matching the requested tolerance (Figure \ref{fig:error_ex1}b). All the solution fields converge super-algebraically with respect to the problem size
(Figure \ref{fig:error_ex1}c).

\begin{figure}[t]  
\centering
\begin{tabular}{ccc}
(a) & (b) & (c)\\
\includegraphics[width=0.29\textwidth]{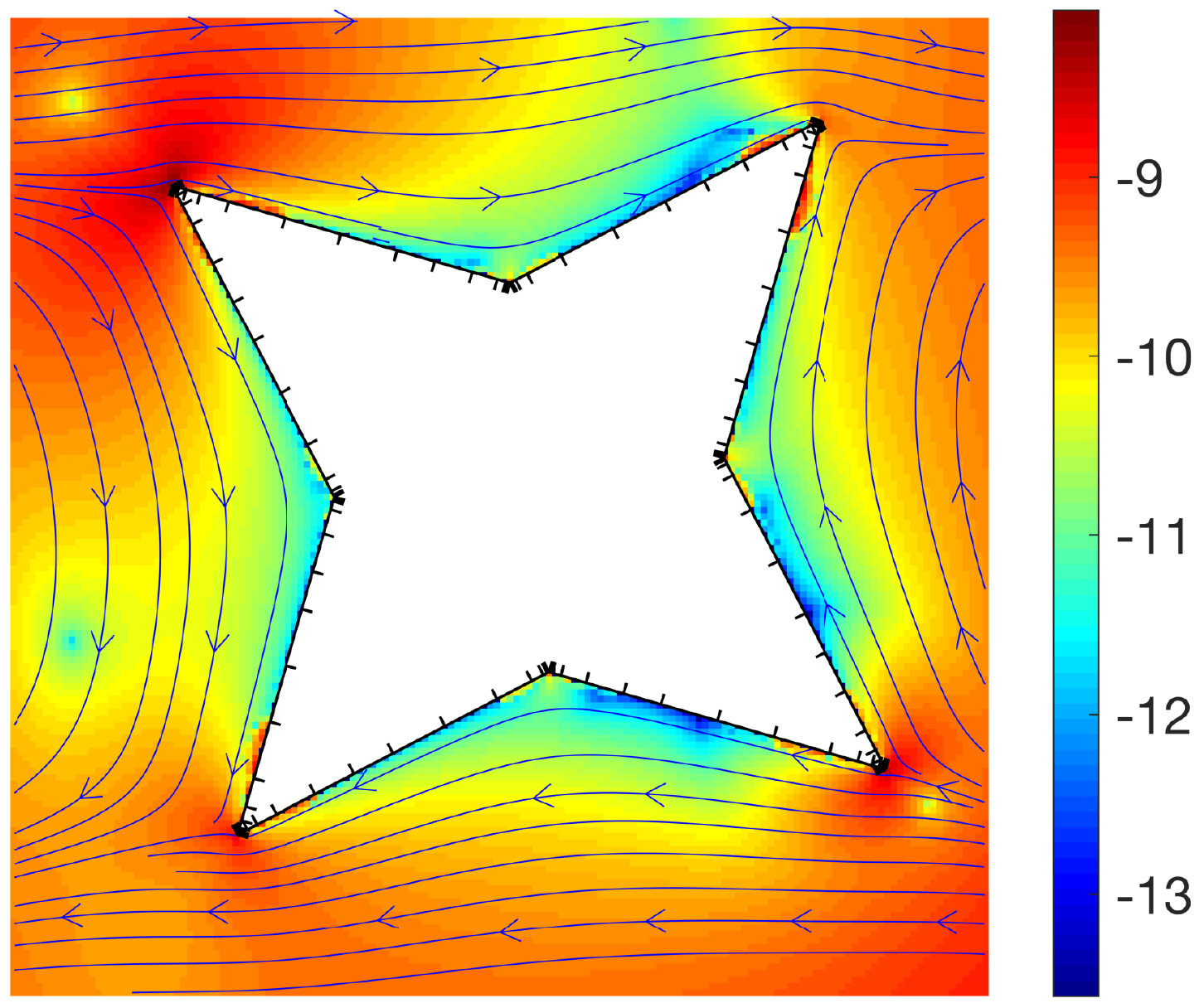}  & \includegraphics[width=0.33\textwidth]{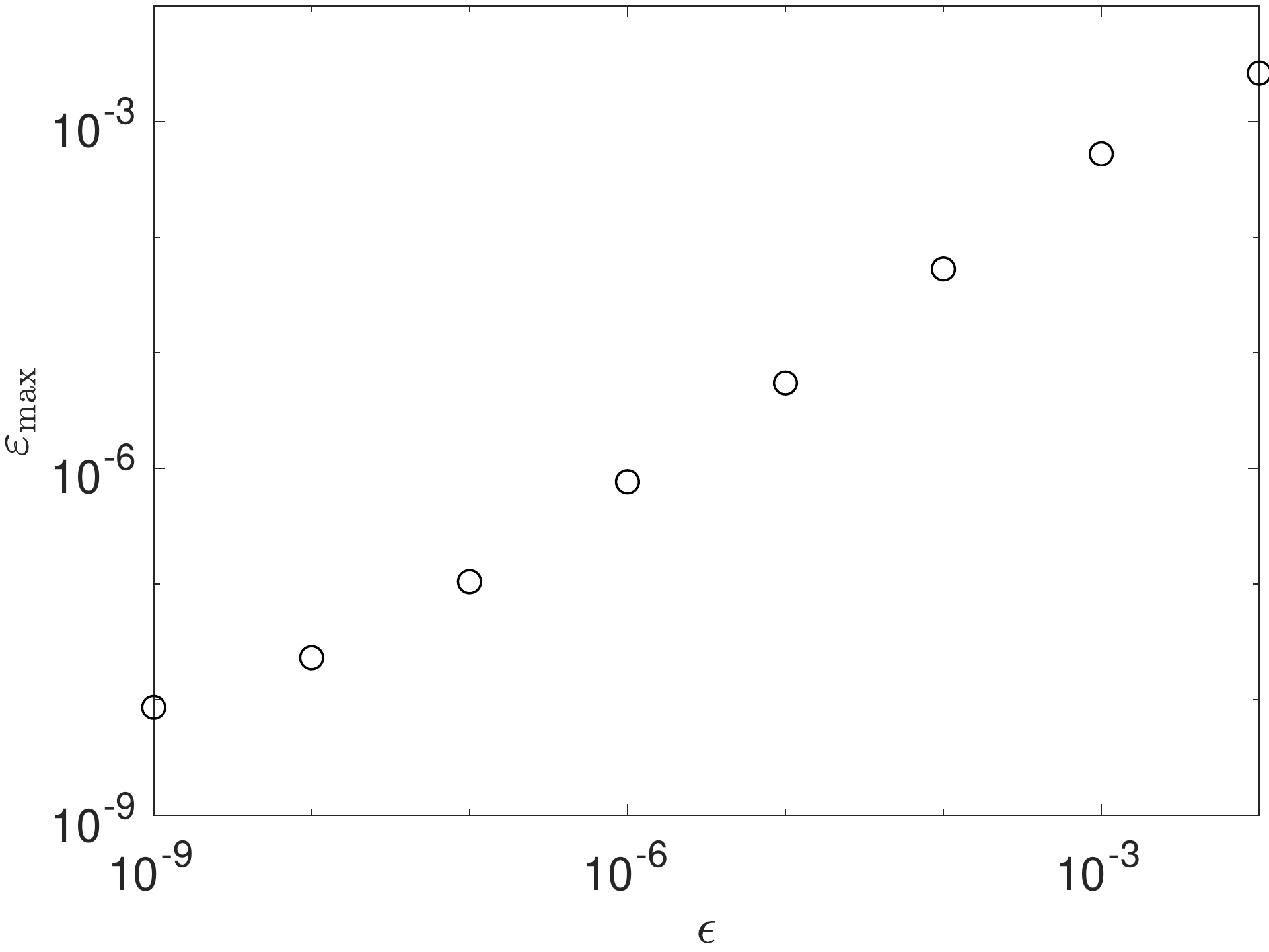} & \includegraphics[width=0.33\textwidth]{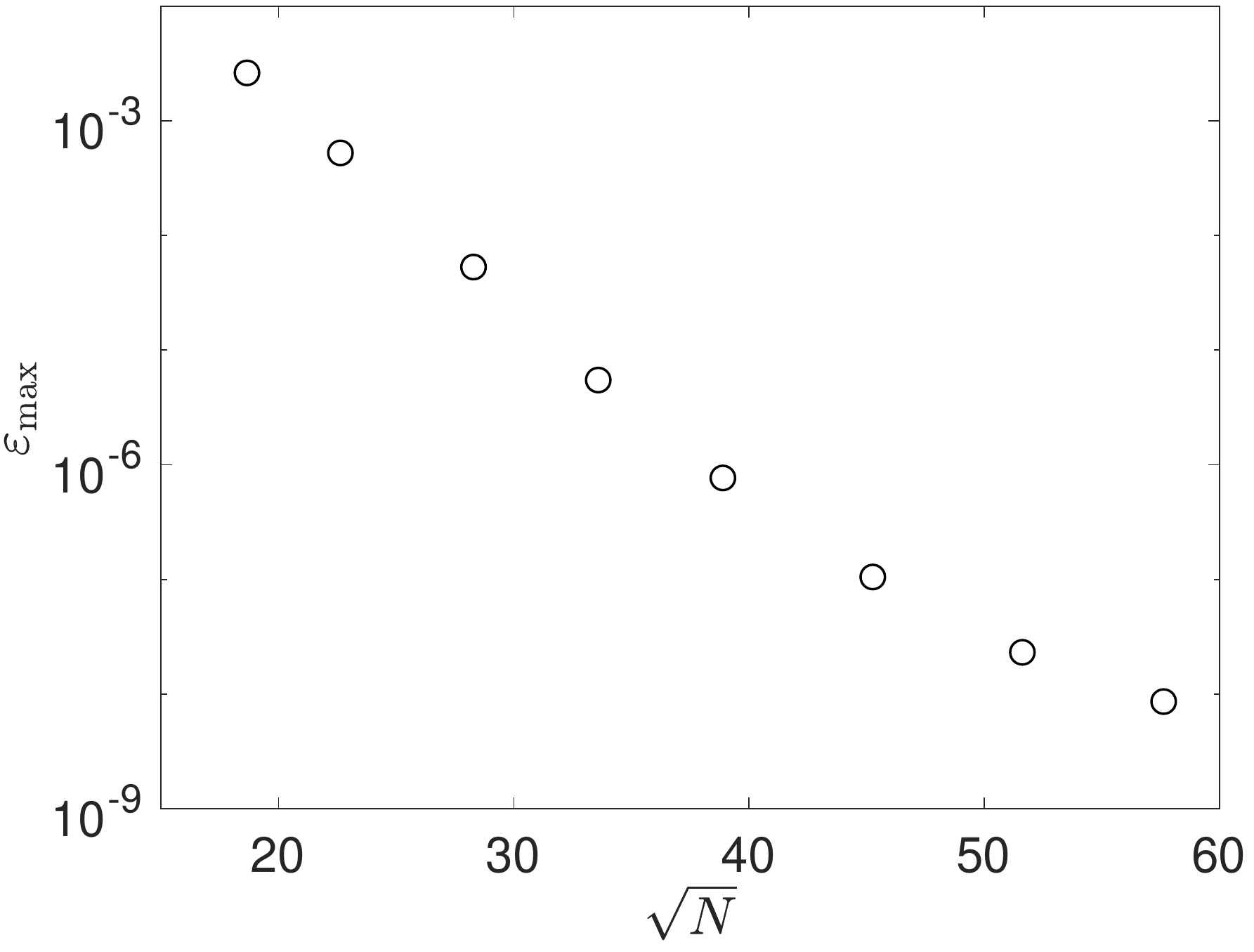}
\end{tabular}
\caption{(a) Linear shear flow around a shuriken-shaped island with 8 corners. Streamlines of the flow and panel endpoints (small segments) are shown. Color represents the $\log_{10}$ error of the velocity computed under tolerance $\epsilon=10^{-10}$ (resulted in  $2N=6640$ degrees of freedom).
(b) Convergence of the maximum error $\varepsilon_\mathrm{max}$ versus requested tolerance $\epsilon$.
(c) Convergence of the maximum error $\varepsilon_\mathrm{max}$ versus the square root of the number of nodes; root-exponential convergence would result in a straight line.
\label{fig:error_ex2}}
\end{figure}

\paragraph{Example 2. Domain with corners.} The smooth geometry in Example 1 is now replaced with a non-convex polygon. Figure \ref{fig:error_ex2} shows a linear shear flow around a ``shuriken'' domain with eight corners, the outer four of which are reentrant (with respect to $\Omega$) corners of angle $\theta=1.74\pi$. With $\alpha = 0.5$ for the flatter corners and $\alpha=1.1$ for the sharper ones, our scheme achieved accuracies that are matching the requested tolerance (Figure \ref{fig:error_ex2}b). Note that the convergence with respect to problem size is super-algebraic  (Figure \ref{fig:error_ex2}c), and consistent with root-exponential convergence, as expected for problems with corner singularities
(see Remark~\ref{r:rootexp}). We used the numerical solution obtained using  $\epsilon=10^{-10}$ as the reference solution.

\begin{figure}[t]  
\centering
\begin{tabular}{ccc}
(a) & (b) & (c)\\
\includegraphics[width=0.3\textwidth]{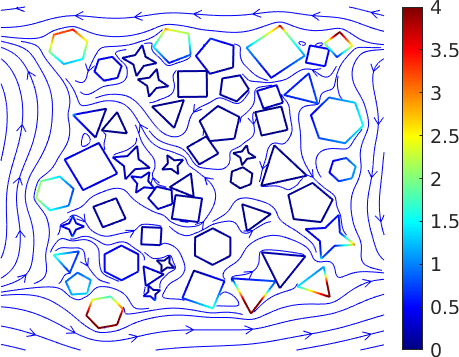}  & \includegraphics[width=0.32\textwidth]{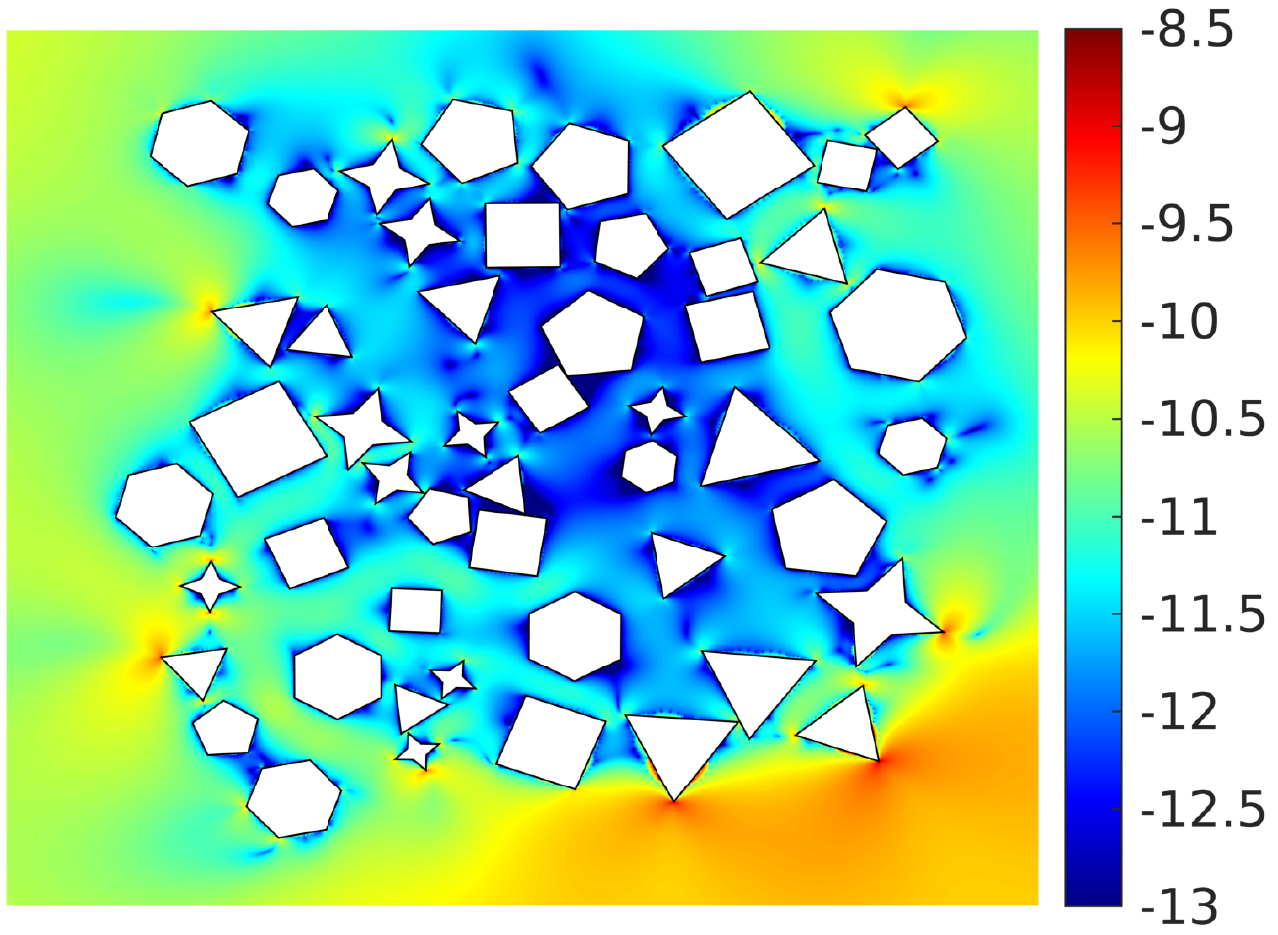} & \includegraphics[width=0.32\textwidth]{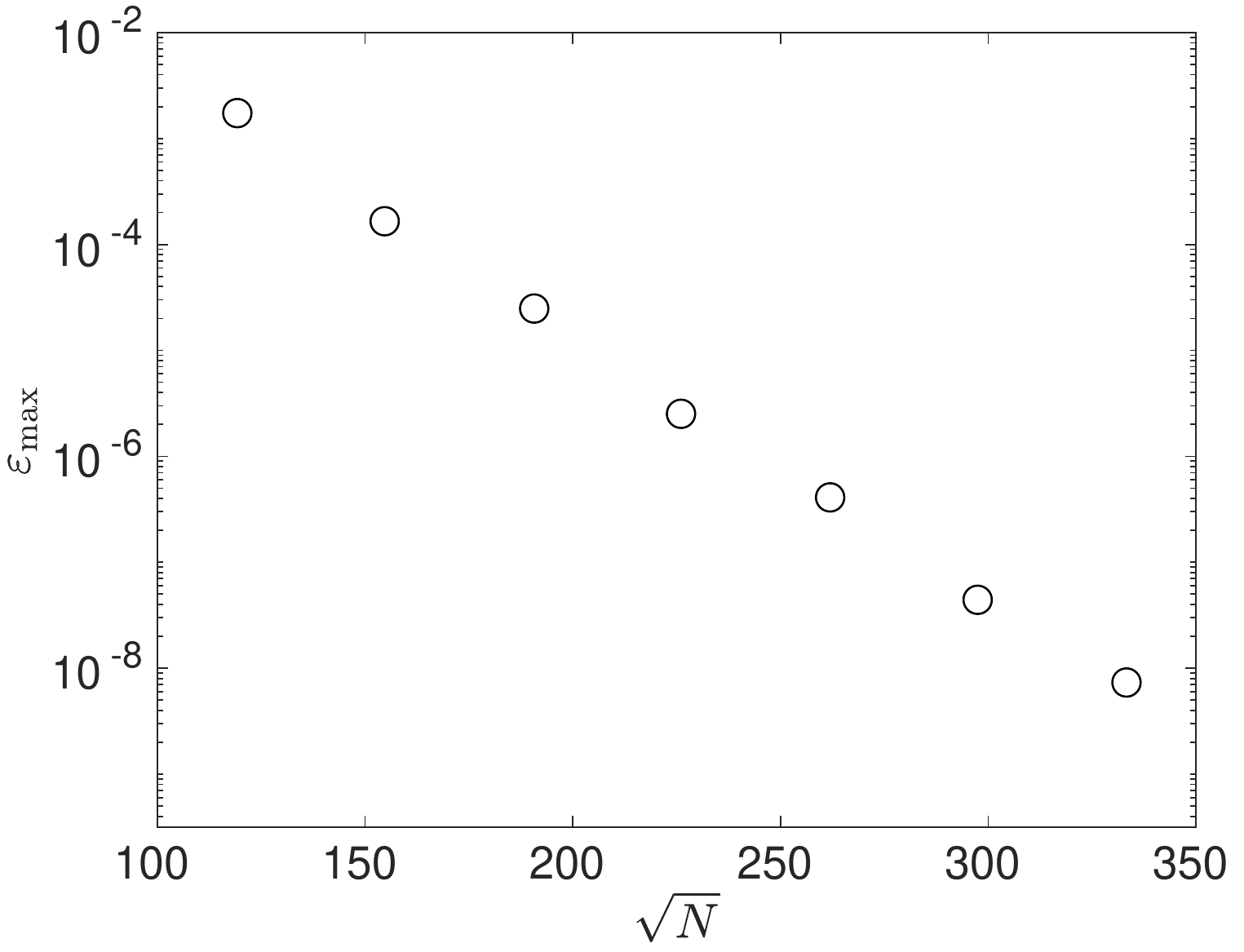}
\end{tabular}
\caption{(a) Streamlines of a shear flow past $50$ randomly generated polygonal islands with a total number of 253 corners. Color on the polygon boundaries indicate the magnitude of density $\sigma$. 
(b) $\log_{10}$ of absolute error of the velocity, computed using $2N = 222140$ degrees of freedom. Error is measured on a $1000\times 1000$ grid (spacing $\Delta x \approx 8.3\times10^{-3}$) by comparing to the solution obtained with $\epsilon=10^{-10}$.
(c) Convergence of the maximum error $\varepsilon_\mathrm{max}$ versus the square root of the number of nodes.
\label{fig:error_ex3}}
\end{figure}

\paragraph{Example 3. Multiple polygonal islands.} This example models a porous media flow through a collection of non-smooth, non-convex and closely packed boundaries: we set up $50$ polygonal islands with a total number of 253 corners (Figure \ref{fig:error_ex3}). The computational domain has width $\approx8$ and the closest distance between the polygons is about $10^{-2}$. The background flow
is the same as in the previous examples. With $\alpha=0.75$ and $\lambda=2$ for all corners, the convergence (Figure~\ref{fig:error_ex3}c) is similar to the single-polygon island example (Figure~\ref{fig:error_ex2}b--c), achieving more than 8 digits using approximately $800$ degrees of freedom per corner. This demonstrates the robustness of our adaptive scheme, that is, the performance is as good for a more complex example as for a simple one.

\paragraph{Example 4. Artificial vascular network.} We now turn to the example shown in Figure \ref{fig:vasnet}. We construct an artificial vascular network (with 378 corners) that mimics those observed in an eye of a zebra fish \cite{alvarez2007genetic}. 
The flow in this network is driven by a uniform influx from the circular wall at the center and a uniform outflux at the outer circular wall, such that the overall volume is conserved; all other boundaries have a no-slip condition. 
We solve the BIE for this problem using GMRES with a block diagonal preconditioner consisting of the diagonal panel-wise blocks of the BIE system itself (i.e., the self-evaluation blocks for each panel).
The FMM is used for applying the matrix and for final flow evaluations;
see Remark~\ref{r:fmm}.
The sparse correction matrix (see Section~\ref{s:close_corr})
is applied via MATLAB's single-threaded built-in matrix-vector
multiplication;
its rows have been precomputed as described in Section~\ref{s:matel}.
All computations are done on an $8$-core $4.0$ GHz Intel Core i7 desktop.

\begin{figure}[t]
\centering
\includegraphics[width=\textwidth]{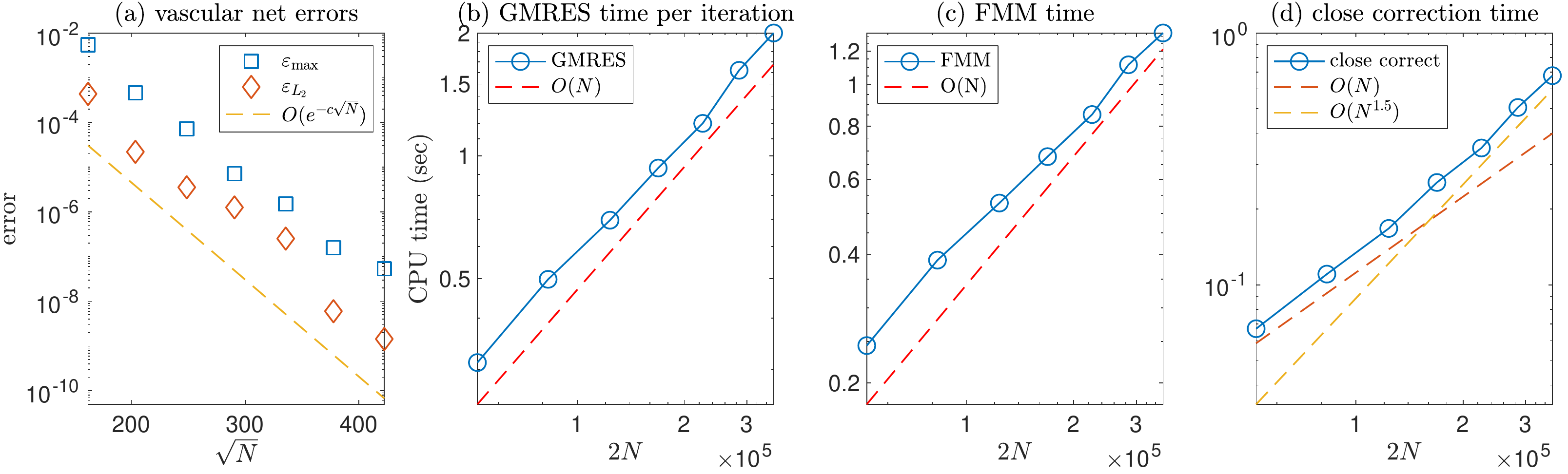}
\caption{Error and timing of solving the vascular network BVP. (a) Convergence of errors in log-linear scale. (b) Log-log scale plot of the total CPU time per GMRES iteration, which consists of the FMM time, shown in (c), and the close correction time, shown in (d). }\label{fig:vasnet_err_timing}
\end{figure}

Table \ref{tab:vasnet} shows, for various tolerance $\epsilon$, the relative $L_2$-error $\varepsilon_{L_2}$, the relative maximum error in velocity $\varepsilon_{\mathrm{max}}$, the total number of panels used $n_\Lambda$,
the number of degrees of freedom $2N$, memory (RAM) usage, the number of GMRES iterations and time used, the setup time for precomputing the close-correction matrices, and the percentage time for applying Stokes FMM  during GMRES. Several observations are in order:
\begin{enumerate}[1)]
\item Both $\varepsilon_{L_2}$ and $\varepsilon_\mathrm{max}$ decay super-algebraically with the number of degrees of freedom; this data is plotted in Figure~\ref{fig:vasnet_err_timing}a.
The closeness between $\varepsilon_{L_2}$ and $\varepsilon_\mathrm{max}$ shows that our scheme has achieved high accuracies near the sharper reentrant corners (hard) that are similar to those near the smooth edges (easy). 
This error analysis remains valid even in the zoomed-in high-resolution error plots in Figure~\ref{fig:vasnet}.
Furthermore, the convergence performance of this example is the same as the previous two examples---we achieved more than 8 digits, with a ratio $\frac{\text{degrees of freedom}}{\text{\#corners}}\approx943$, which is similar to the ratios in Example 2 (830) and Example 3 (878). This once again demonstrates the robustness of our overall scheme to problem complexity.
\item  The number of GMRES iterations increases only because we have requested smaller tolerance. Each additional digit needs about $100$ more iterations.
The GMRES convergence {\em rate} is stable, which demonstrates that our second kind BIE formulation is well-conditioned even in the presence of corners.
\item The fact that the Stokes FMM time is the main cost shows that our algorithm has achieved close to optimal efficiency. 
The slight decrease of the percentage FMM times at smaller $\epsilon$ is due to the fact that the FMM time grows only linearly with $N$, while the close evaluation matrix-vector multiplication time grows like $O(N^{3/2})$. The latter estimate is obtained as follows. The number of matrix-vector multiplications grows like $O(n_\Lambda)=O(\log 1/\epsilon)=O(p)$, where each matrix-vector product takes $O(p^2)$ time. Note that $N=n_\Lambda\times p = O(p^2)$, so the total close evaluation time grows as $O(p^3)=O(N^{3/2})$. (See Figure~\ref{fig:vasnet_err_timing}b--d.)
\end{enumerate}

\begin{table}[t]
  \centering
	\begin{tabular}{c|p{1.5cm}|p{1.5cm}|p{1cm}|p{1cm}|p{1.2cm}|p{1.31cm}|p{1.3cm}|p{1cm}|p{1cm}}
	\hline
	$\epsilon$ & $\varepsilon_{L_2}$ & $\varepsilon_\mathrm{max}$ & $n_\Lambda$ & $2N$ & RAM used(gb) & GMRES iteration & GMRES time(s)& setup time(s) & \% FMM\\
	\hline            
	1e-03 &  4.34e-04 & 5.43e-03  & 6549  & 52392  & 2.3  & 796 &  248 & 48 & 78.50\\
	1e-04 &  2.20e-05 & 4.58e-04  & 8281 &  82810 & 2.9 & 919 & 458 & 63 & 77.86 \\
	1e-05 &  3.55e-06 & 7.21e-05  & 10301  & 123612  & 3.7  & 1091 & 759 & 84 & 75.92\\
	1e-06 &  1.26e-06 & 7.15e-06  & 12061 & 168854  & 5.0  & 1282 & 1197 & 106 & 72.69\\
	1e-07 &  2.53e-07 &  1.51e-06 & 14079 & 225264  & 6.9 & 1390 & 1670 & 135 &  70.91\\
	1e-08 &  6.01e-09 &  1.58e-07 & 15839 & 285102 & 9.0 & 1501 & 2433 & 164 & 68.77\\
	1e-09 &  1.44e-09 &  5.34e-08 & 17829 & 356580 & 12 & 1597 & 3195 & 204 & 66.16 \\
	\hline
	\end{tabular}
	\caption{\em Results and statistics of solving the BVP in the vascular network in Figure \ref{fig:vasnet} for various tolerance $\epsilon$. Errors $\varepsilon_\mathrm{max}$ and $\varepsilon_{L_2}$ are measured on a $2160\times 2160$ grid (spacing $\approx 2.5\times10^{-3}$) by comparing to the solution obtained at $\epsilon=10^{-10}$. CPU time and RAM used are measured using \cite{memorygraph}. }\label{tab:vasnet}
\end{table}

\begin{figure}[t]
\centering
\includegraphics[height=1.6in]{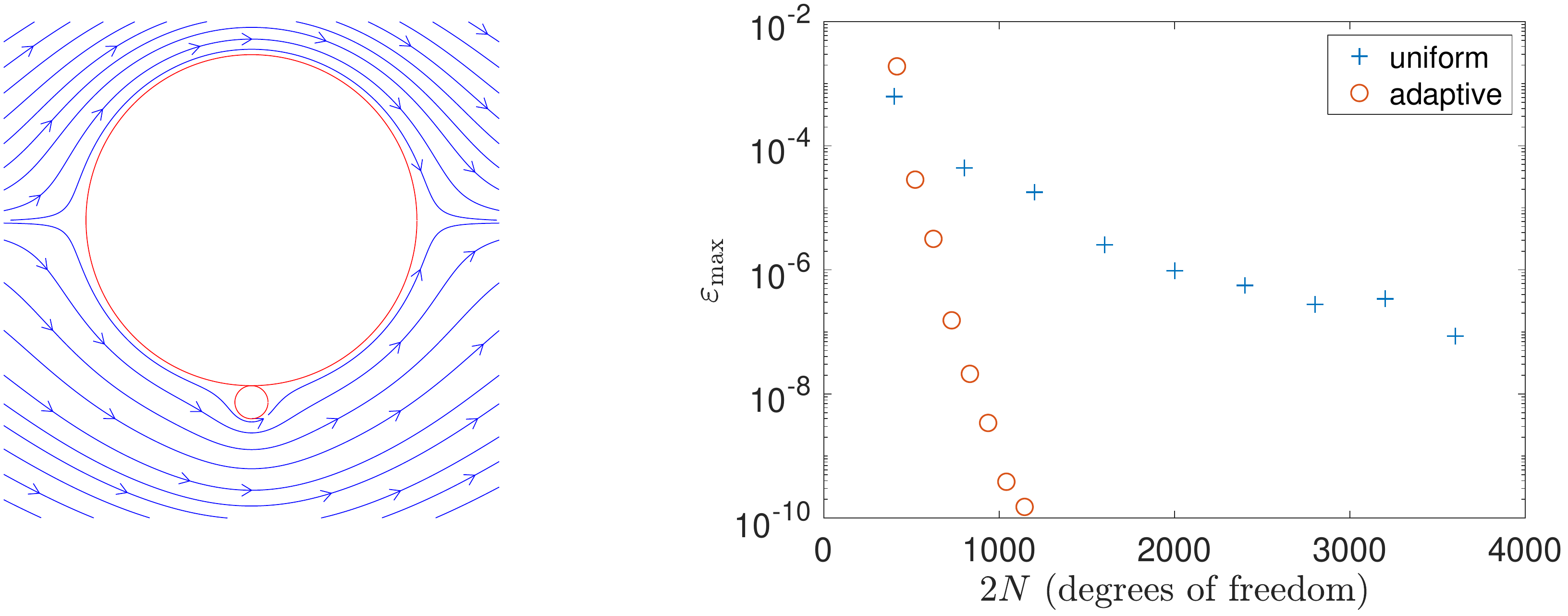}
\caption{Convergence of a uniform flow past two touching disks that are $d=10^{-6}$ apart, and whose radii are 1 and 0.1. The required number of unknowns in the adaptive scheme is much less than the global scheme with uniform resolution; see Example 5.}\label{fig:globalVpanel}
\end{figure}

\paragraph{Example 5. Uniform versus adaptive for close-to-touching curves.} Finally, to demonstrate the advantage of using an adaptive scheme over a uniform discretization, let us consider a uniform flow past two close-to-touching disks (Figure \ref{fig:globalVpanel}).
The background flow is a constant $u_\infty = (1,0)$, the separation is $d=10^{-6}$, and the radii 1 and 0.1.
For the uniform-resolution scheme we use a global periodic trapezoid grid on each circle, where, in order to have similar node spacings, the larger circle has $9$ times as many points as the smaller one.
Here, global close-evaluation is done using the spectrally accurate quadrature from \cite{lsc2d}.
On the other hand, the adaptive quadrature uses a grid that is determined by our adaptive refinement scheme of Section~\ref{sc:amr}, with one modification
that improves the scaling in the number of refined panels (see Remark~\ref{r:rootdist_rule}).
We observe that, for more than 4 accurate digits, the number of unknowns required by the adaptive scheme is much less than that of the uniform-resolution scheme (Figure \ref{fig:globalVpanel}).
The smoothness of the density function discussed in the remark below
suggests that, at fixed $\epsilon$,
the uniform scheme (and also the original refinement scheme) needs
$N=O(1/\sqrt{d})$ unknowns,
whereas the modified adaptive scheme needs only $N=O(\log(1/d))$.
The latter is close to optimal, and is what we recommend for closely-interacting
curves. (See also Examples 4 and 6 in \cite{barnett-ls2p2d} for a ``globally adaptive'' variant.)

\begin{remark}[Refinement at close-to-touching smooth surfaces]\label{r:rootdist_rule}
For viscous flow in the region between two smooth curves
separated by a small distance $d$,
asymptotic analysis gives that the width of the ``bump''
in fluid force scales as $O(\sqrt{d})$ \cite{sangani1994inclusion}.
By dimensional analysis, if the sum of the two curvatures of the
surfaces near the contact point is $\kappa$, then
the width in fact scales as $O(\sqrt{d/\kappa})$.
Assuming that this also applies to the density $\sigma$,
this suggests a looser criterion for refinement:
panels should be refined only when they are longer than this width scale.
This allows panels to come much closer than their length, without being
refined. In the case of close smooth curves,
the test in line 13 of Algorithm~\ref{alg:glob_refine} can thus be modified to
$(c'\sqrt{\kappa}S)^2>d$. We find that the constant $c'=0.7$ achieves
the requested tolerance.
The resulting $N$ can be estimated as follows.
Setting $\kappa=1$ for simplicity,
a generic local model of the separation is $h(t) \approx d +t^2$ as a function
of parameter $t$, and $n_\Lambda = O\bigl(\int_{-1}^1 dt/S(h(t))\,\bigr)$
where $S(h)$ is the local panel size as a function of separation.
The original refinement scheme, $S(h) = O(h)$, thus gives
$n_\Lambda = O(d^{-1/2})$, whereas the modified $S(h)=O(\sqrt{h})$
gives $n_\Lambda = O(\log(1/d))$.
\end{remark}

\section{Conclusions}\label{sc:conclusions}
We have presented a set of panel quadrature rules for accurate evaluation of single- and double-layer Stokes potentials and their associated pressure and traction boundary integrals. They can be used for targets that are either on or off the boundary, and can be located arbitrarily close to it. In addition, we formulated an adaptive panel refinement procedure that sets the length of panels on the boundary (``$h$-adaptivity'') and the overall degree of approximation $p$ (``$p$-adaptivity'') required to achieve a user-prescribed tolerance. We demonstrated via numerical experiments that our algorithm achieves super-algebraic convergence even for complex geometries with corners, and that the CPU time grows linearly with problem size, and is dominated by the cost of FMMs for large-scale problems. More sophisticated quadratures and techniques designed for corner singularities, such as the RCIP \cite{helsing2013solving, helsing2018integral} (or the work of \cite{rachh2017solution}), are expected to further improve the performance of our BIE solver.
It is also expected that adapting $p$ on a per-panel basis
(i.e.\ full $hp$-adaptivity)
would reduce the total number of degrees of freedom needed,
although only by a factor less than two.
Applications of our work include providing design tools for rapid prototyping of microfluidic chips (for cell sorting, mixing or other manipulations e.g., \cite{kabacaouglu2018optimal}), shape optimization (e.g., \cite{bonnet2019shape}) and  simulating cellular-level blood flow in microvasculature.

We envision building a fast 2D particulate flow software library by utilizing the algorithms developed in this work for the fixed complex geometries (such as microfluidic chips or vascular networks) and our global close evaluation schemes developed in \cite{lsc2d} for the moving rigid or deformable particles (such as colloids, drops or vesicles). Another key ingredient would be a {\em fast direct solver} for solving the BIEs on the fixed geometries, similar to that developed in \cite{marple-periodic}, wherein the boundary integral operators were compressed by exploiting their low-rank structures, inverted as a precomputation step, and applied at an optimal $\bigO(N)$ cost at every time-step of the particulate flow simulation as particles move through the fixed geometry. One open research question in this context is: {\em Can we update the compressed representations as the boundary panels are refined (or coarsened) without rebuilding them?} A similar question was recently investigated in \cite{zhang2018fast}, where the
authors report a $3\times$ speedup when locally perturbing the
geometry. We plan to explore their approach and report its performance in the context of our adaptive panel refinement procedure.

\section{Acknowledgements}
We thank Leslie Greengard and Manas Rachh for many useful discussions pertaining to this work. We thank Manas Rachh for his biharmonic FMM code, and Jun Wang for providing Stokes MATLAB wrappers to this code. We acknowledge support from NSF under grants DMS-1719834 and DMS-1454010. The work of SV was also supported by the Flatiron Institute, a division of the Simons Foundation.


\appendix
\begin{appendices}
  \section{Pressure and traction in terms of contour integrals}
  \label{a:trac}

Here, we give formulae for the traction vector $\bf T$ induced at a target point
with given surface normal, and the
associated pressure field $p$, when the velocity field is represented by a Stokes single or double layer potential.
The goal is to write the traction and pressure in terms of the four contour
integrals of Sec.~\ref{s:cont}, to which close-evaluation methods of
Sec.~\ref{s:eval} may then be applied.
This enables uniformly accurate force calculations on bodies, or
solution of traction BVPs.
We use the notation of Sec.~\ref{sc:problem}:
recall that $r = x - y$, $n_x$ and $n_y$ are the normal vectors at the target $x$ and source $y$ respectively, $\rho=|r|$, and $I$ denotes the $2\times 2$ identity operator.

We first consider the single layer potential \eqref{sto_slp}.
Its traction is
\begin{equation}
(\mathbf{T}^S\sigma)(x) = -\dfrac{1}{\pi}\int_{\Gamma}\dfrac{r\cdot n_x}{\rho^2}\,\dfrac{r\otimes r}{\rho^2}\sigma(y)\,ds_y,
\end{equation}
which turns out to be the negative of the Stokes DLP
\eqref{sto_dlp} with $n_y$ replaced by $n_x$.
While \eqref{identity_for_dlp_split} is no longer useful in this case,
we can instead write the traction as
\begin{equation}
(\mathbf{T}^S\sigma)(x) = -\dfrac{1}{\pi}\int_{\Gamma}(r\cdot\sigma(y))\dfrac{r}{\rho^4}(r\cdot n_x)\,ds_y
\end{equation}
and use the slightly different identity
\begin{equation}
\nabla_x\left(\dfrac{r\cdot \sigma}{\rho^2}\right) = \dfrac{\sigma}{\rho^2} - (r\cdot \sigma)\dfrac{2r}{\rho^4}
\end{equation}
to write the traction kernel as
\begin{equation}
\begin{aligned}
(\mathbf{T}^S \sigma)(x)\;
=\;& {-1\over2\pi}\int_{\Gamma}{  r\cdot n_{x}\over\rho^2}\sigma\,ds_{y}  + {1\over2\pi}(x\cdot n_x)\Del\int_{\Gamma}{ r\cdot\sigma \over\rho^2}\,ds_{y} \\
 & - {1\over2\pi}n_{x,1}\Del\int_{\Gamma}{ r\cdot\sigma \over\rho^2}y_1\,ds_{y} - {1\over2\pi}n_{x,2}\Del\int_{\Gamma}{ r\cdot\sigma \over\rho^2}y_2\,ds_{y}
~,
\end{aligned}
\label{slp_traction_via_lps}
\end{equation}
where $(n_{x,1}, n_{x,2})=:n_x$ are the two components of $n_x$. As did in the case of velocity potentials, we can concisely write \eqref{slp_traction_via_lps} as
\begin{equation}
  \mathbf{T}^S\sigma = \big((\mathcal{S}\sigma_1)_n+i (\mathcal{S}\sigma_2)_n\big)  + \frac{1}{2\pi}\left(\text{Re}(x/n_x)\overline{I_H(\sigma/n_y)} - n_{x,1}\overline{I_H(\sigma y_1/n_y)} - n_{x,2}\overline{I_H(\sigma y_2/n_y)}\right),
  \label{tracIH}
\end{equation}
where all the $\mathbb{R}^2$ vectors are now understood as complex numbers in $\mathbb{C}$, the over line in $\overline{I_H(\cdot)}$ denotes the complex conjugate of $I_H(\cdot)$ and the dot product $x\cdot n_x = \text{Re}(x/n_x)$ is due to the fact that $1 = |n_x|^2 = n_x\overline{n_x}$.

The single layer pressure associated to \eqref{sto_slp} is
\begin{equation}
(P^S\sigma)(x) = \dfrac{1}{2\pi}\int_{\Gamma}\dfrac{r\cdot \sigma}{\rho^2}\,ds_y~,
\end{equation}
which again in the complex plane can be written as 
\begin{equation}
  P^S\sigma = \text{Re}\frac{i}{2\pi}I_C(\sigma/n_y).
  \label{presIC}
\end{equation}

We now turn to the Stokes double layer potential \eqref{sto_dlp}.
The traction kernel and its associated pressure kernel are given by \cite[(5.27)]{liu2009fast} \cite[(3.37)]{barnett-ls2p2d},
\begin{equation}
\setlength{\jot}{10pt}
    \begin{aligned}
        (\mathbf{T}^D\sigma)(x) =& \frac{\mu}{\pi} \int_{\Gamma} \left( -8\frac{r\otimes r }{\rho^6} (r\cdot n_x) (r\cdot n_y)
        +\frac{r\otimes n_x }{\rho^4}(r\cdot n_y) + \frac{r\otimes r}{\rho^4}(n_x\cdot n_y) \right.\\
        &\left. \hspace{0.9in} +I\frac{1}{\rho^4}(r\cdot n_x) (r\cdot n_y) + \frac{n_y\otimes r}{\rho^4}(r\cdot n_x) +\frac{n_x\otimes n_y}{\rho^2} \right)\sigma(y)\,ds_y \\
       (P^D\sigma)(x) =& \frac{\mu}{\pi}\int_{\Gamma}\left( -\frac{n_y\cdot \sigma(y)}{\rho^2} + 2\frac{r\cdot \sigma(y)}{\rho^4}\left(r\cdot n_y\right)\right)\,ds_y.
    \end{aligned}
\end{equation}
The corresponding boundary integral operators  $(\mathbf{T}^D \sigma)(x)$ and $(P^D \sigma)(x)$ are hyper-singular. We can easily derive the following equation expressing this operator in terms of the Laplace double layer potential:
\begin{equation}
\setlength{\jot}{10pt}
    \begin{aligned}
        \frac{1}{\mu}\mathbf{T}^D \sigma = &
            -2\left(x\cdot n_x\nabla\nabla \mathcal{D}[\sigma] - n_{x,1}\nabla\nabla\mathcal{D}[y_1 \sigma] -n_{x,2}\nabla\nabla[y_2\sigma]\right) \\
            &+3I(n_x\cdot \nabla \mathcal{D})[\sigma]  -(n_x\otimes\nabla \mathcal{D})[\sigma]-(\nabla \mathcal{D} \otimes n_x)[\sigma] \\
            & +\begin{bmatrix}1 & \\ & -1\end{bmatrix}
            \left(n_x\cdot \nabla \mathcal{D} + n_x\otimes \nabla\mathcal{D} - \nabla\mathcal{D}\otimes n_x\right)\left[\frac{\bar{n}_y}{n_y}\sigma\right],\\
        P^D \sigma = & -2\mu\left(\frac{\partial}{\partial x_1} \mathcal{D}[\sigma_1] + \frac{\partial}{\partial x_2} \mathcal{D}[\sigma_2]\right)\\
        			= & \frac{\mu}{\pi}\left(\text{Im}\left(I_H(\sigma_1)\right)+\text{Re}\left(I_H(\sigma_2)\right)\right).
    \end{aligned}
\end{equation}

Here $\frac{\bar{n}_y}{n_y}\sigma = [\frac{\bar{n}_y}{n_y}\sigma_1,\frac{\bar{n}_y}{n_y}\sigma_2]^T$, $\nabla$ is short for $\nabla_x$, and $\nabla\nabla$
is the Hessian tensor.
The gradients of Laplace double-layer potentials needed above
are expressed in terms of Hadamard integrals using \eqref{gradD}.
The Hessians are given in terms of supersingular integrals as
follows:
\begin{equation}
\nabla\nabla\mathcal{D}[\sigma](x) \;=\;
\begin{bmatrix}
  \text{Re}\frac{i}{\pi} \left(I_S(\sigma_1)\right)(x) &
  -\text{Im}\frac{i}{\pi} \left(I_S(\sigma_2)\right)(x) \\
  -\text{Im}\frac{i}{\pi} \left(I_S(\sigma_1)\right)(x) &
  -\text{Re}\frac{i}{\pi} \left(I_S(\sigma_2)\right)(x)
\end{bmatrix}.
\end{equation}
The close evaluation formulae for these are in Sec.~\ref{sec:IH}.

To validate the above formulae, we include in Fig.~\ref{fig:error_ex1}(b--c)
the convergence of the maximum error in pressure and traction for the smooth domain of Example 1 from Sec.~\ref{sc:results}.
The convergence rate is very similar to that of velocity 
albeit a loss of 1--2 digits, which is expected due to the extra derivatives.

\end{appendices}

\bibliographystyle{unsrt}
\bibliography{bobi,3dVesicles}

\end{document}